\documentclass{nm}
\renewcommand\thisnumber{x}
\renewcommand\thisyear {2023}
\renewcommand\thismonth{x}
\renewcommand\thisvolume{xx}
\renewcommand\doi{10.4208/nmtma.  }

\usepackage{charter,rotating,enumerate,multirow,algorithm,algorithmic,}
 
\def\nnz{{\rm{nnz}}}
\def\diag{{\rm{diag}}}
\def\PPS{{\rm{PPS}}}

\def\SPPS{{\rm{SPPS}}}

\def\nul {\text{\textsf{null}}}
\def\rank {\text{\textsf {rank}}}
\def\ev {\text{\textsf {eigvec}}}

\newcommand{\norm}[1]{\left\|{#1}\right\|}

\begin{document}

\markboth{    M. Masoudi and D. K. Salkuyeh}{   PPS  iteration methods for solving nonsingular   positive	semidefinite  systems    }
\title{  Positive semidefinite  and positive semidefinite splitting  iteration methods for solving nonsingular   positive	semidefinite  systems     }

\author[1 ]{ Mohsen Masoudi\affil{1}\comma\affil{2},\corrauth and Davod Khojasteh Salkuyeh\affil{3}}
\address{
\affilnum{1}\   Minab Higher Education Center,  University of Hormozgan,  P.O. Box  3995, Bandar Abbas,  Iran\\ 
\affilnum{2}\ Department of Mathematics, University of Hormozgan, P.O. Box 3995, Bandarabbas, Iran \\
\affilnum{3}\ Faculty of Mathematical Sciences, University of Guilan, P.O. Box 1914, Rasht, Iran\\
[2ex]
\rm Received   ; Accepted (in revised version)  }

\emails{
{\tt  m.masoudi@hormozgan.ac.ir}  (M. Masoudi),
{\tt  khojasteh@guilan.ac.ir}  (D. K. Salkuyeh)
}

\begin{abstract}
This article introduces an iterative method for solving nonsingular  positive semidefinite systems of 	linear equations. To construct the iteration process,  the coefficient matrix is split into two   positive semidefinite matrices along with an arbitrary Hermitian positive definite shift matrix.  Several conditions are established to guarantee the convergence of  method and  suggestions are provided for selecting the matrices involved in the desired splitting.	We explore  selection process of the shift matrix and determine the optimal parameter in a specific scenario. The proposed method aims to generalize previous approaches and improve the conditions for  convergence theorems. In addition, we examine  two special cases of this method and compare the induced preconditioners with some state-of-art preconditioners. Numerical examples are given to demonstrate effectiveness of  the presented preconditioners.
\end{abstract}

 
 
\maketitle

\section{Introduction}
Consider  the solution of systems of linear equations  in the form 
\begin{equation}\label{system}
	\mathcal{A} {u}= {b},
\end{equation}
where  ${b}\in \mathbb{C}^{n}$ and $\mathcal{A}\in \mathbb{C}^{n \times n}$ is large, sparse, nonsingular  positive semidefinite (PSD), i.e.,   $\mathcal{A}+ \mathcal{A^*}$   is Hermitian positive semidefinite (HPSD).	{
	Systems of this type arise in numerous scientific and engineering applications, such as computational fluid dynamics and constrained optimization \cite{bai752006,benzi142005,rees322010}.}
Various  methods can be utilized to solve a system of linear equations when $\mathcal{A}$   is assumed to be positive definite (PD),  i.e.,   $\mathcal{A}+ \mathcal{A^*}$   is Hermitian positive  definite (HPD).
For instance, the Hermitian and skew-Hermitian splitting (HSS) method was introduced by Bai et al. in \cite{bai242003}. 
This method uses the Hermitian and skew-Hermitian splitting of the matrix $\mathcal A$ as $\mathcal A=\mathcal H+\mathcal S$,  where $\mathcal H=\frac{1}{2} (\mathcal A+\mathcal A^* )$ and $\mathcal S=\frac{1}{2} (\mathcal A-\mathcal A^*) $. The iteration process is then based on this splitting:
\begin{equation*}
	\begin{cases}
		(\alpha \mathcal I+ {\mathcal H}  ) {u}^{k+\frac{1}{2}}= (\alpha \mathcal I -  {\mathcal S}  )  {u}^{k}+{b},\\
		(\alpha \mathcal I + {\mathcal S}  ) {u}^{k+1}= ( \alpha \mathcal I - {\mathcal H
		}  )  {u}^{k+\frac{1}{2}}+{b},	
	\end{cases}
\end{equation*}
in which $\alpha>0$,  $\mathcal I$ is identity matrix and    $ u^0 $ is an arbitrary initial guess. 

It is shown in \cite{bai242003} that the HSS method converges unconditionally to the unique solution of the system \eqref{system} when A is  positive definite.
Due to the method's efficiency, it has attracted significant attention from researchers, leading to the developing several extensions. The preconditioned version of the HSS method, known as PHSS, was proposed by Bai et al. in \cite{bai982004}.
The modified HSS (MHSS) and its preconditioned version  were introduced by Bai et al. for solving symmetric complex systems \cite{bai872010,bai562011}.	Two generalized versions of HSS (GHSS) were introduced in \cite{benzi312009,zhou2712015}. 
The positive definite and skew-Hermitian splitting (PSS) method, introduced by Bai et al. in \cite{bai262005}, was proposed to solve  positive definite linear systems. It was proved in \cite{bai262005} that the PSS method converges unconditionally to the unique solution of such systems. By further generalizing the HSS iteration method, the normal and skew-Hermitian splitting (NSS) iteration method was proposed  in   \cite{bai142007}. A generalization of PSS  (GPSS)  was introduced by Cao et al.  in \cite{cao22012}. An extension of this method was stated in \cite{masoudi792020}  by Masoudi and Salkuyeh.  Huang and Ma in  \cite{huang342016} presented the  positive definite and semidefinite splitting  method. The shift-splitting (SS) iterative  method and the corresponding induced preconditioner was  presented by Bai in \cite{bai242006}. 

Similar  to the HSS method, we introduce an iterative method, called \textit{{Positive semidefinite and Positive semidefinite Splitting (PPS)}} method, for solving the system \eqref{system}.  The proposed method employs a splitting of the form $\mathcal{A}=  {\mathcal P}_1+{\mathcal P}_2$,  where  ${\mathcal P}_1$ and ${\mathcal P}_2$ are two  positive semidefinite matrices, combined with an arbitrary HPD matrix $\Sigma$ as a shift. Our method encompasses various existing approaches and  enhances their convergence theorems. We establish conditions under which the method is guaranteed to converge. An advantage of  our method over the other methods is the flexibility in selecting matrices  ${\mathcal P}_1$ and ${\mathcal P}_2$. Some ideas concerning choosing the matrices  ${\mathcal P}_1$ and ${\mathcal P}_2$, and the shift matrix $\Sigma$ are discussed.

Throughout this paper, we use the following notations. The set of $ n\times m $ complex matrices is denoted by $ \mathbb{C}^{n\times m} $ and $ \mathbb{C}^n= \mathbb{C}^{n\times1} $ represents the set of 
$n$-dimensional complex vectors.
 {The symbol   $ \mathcal{I} $  denotes the identity matrix.}
For $  A \in \mathbb{C}^{n\times n} $, the notations   $\nul{(A)}$,  $\rank{(A)}$,  $\rho(A)$,  $\sigma(A)$, $\norm{A}$ and  \ev($A$)   stand for the  null space,    rank,  spectral radius, spectrum,   Euclidean  norm   and the set of all nonzero eigenvectors  of the matrix $A$, respectively. If $Q \in \mathbb{C}^{n\times n} $ is  an HPD matrix and $Q=U^*DU$ is  the {spectral decomposition} of  matrix $Q$, then we define $Q^{\frac{1}{2}}=U^*D^{\frac{1}{2}}U$.

The remainder of this paper is organized as follows. In Section \ref{sec:PPS}, we introduce the PPS iterative method and its induced preconditioner.	Our research, outlined in Section \ref{sec:conver},
focuses on analyzing the necessary and sufficient conditions for the convergence of the method, as well as establishing the specific requirements for its existence. The selection of the matrix $\Sigma$ and determination of the optimal parameter are investigated in Section \ref{sec:sigma}.	In Section \ref{sec:p1p2}, we discuss choosing the matrices ${\mathcal P}_1$ and ${\mathcal P}_2$. To demonstrate the theoretical results, we present numerical examples in Section \ref{sec:examp}. Finally, the paper  {concludes with} Section \ref{sec:conc}.

\section{The PPS method}\label{sec:PPS}
Assume that the matrix $\mathcal{A}$  has a splitting of the following form
\begin{align}
	\label{orgsplit}
	\mathcal{A}=   {\mathcal P}_1+   {\mathcal P}_2,
\end{align}
where ${\mathcal P}_1$ and ${\mathcal P}_2$ are  PSD.
In this case, we say that the matrix $\mathcal{A}$ has a \textit{{Positive semidefinite and Positive semidefinite  (PP)}} splitting. 
Consider the following two splittings 
\begin{align*}
	\mathcal{A}	=&  (\Sigma+ {\mathcal P}_1  ) -  (\Sigma-  {\mathcal P}_2  ) \\=&  ( \Sigma+  {\mathcal P}_2 ) -  (\Sigma - {\mathcal P}_1 )	, 
\end{align*}
where {$\Sigma$   is an   arbitrary HPD matrix. }
Obviously, both of the matrices  $\Sigma+ {\mathcal P}_1$  and $\Sigma+ {\mathcal P}_2$ are nonsingular. Now, similar to the   HSS  method \cite{bai242003},	we establish the PPS method as follows:
\begin{align}
	\label{seq}
	\begin{cases}
		(\Sigma+ {\mathcal P}_1  ) {u}^{k+\frac{1}{2}}= (\Sigma -  {\mathcal P}_2  )  {u}^{k}+{b},\\
		(\Sigma + {\mathcal P}_2  ) {u}^{k+1}= ( \Sigma - {\mathcal P}_1  )  {u}^{k+\frac{1}{2}}+{b},	
	\end{cases}
\end{align}
where ${u}^0\in \mathbb{C}^n$ is a given initial guess. 

 {Note that all the PD, HPSD or skew-Hermitian matrices  are   PSD. }
If   ${\mathcal P}_2$ is a skew-Hermitian and $  {\mathcal P}_1$ is an HPD (resp. PD), then  the   PPS method reduces to the    PHSS \cite{bertaccini992005} (resp. Extend of PSS (EPSS) method \cite{masoudi792020}). Now, 	let $\Sigma=\alpha \mathcal I$.  {If $ {\mathcal P}_2=0$}, then the PPS method becomes the shift-splitting  (SS) method \cite{bai242006}. If  $ {\mathcal P}_2$ is a skew-Hermitian and $ {\mathcal P}_1$ is an HPD (resp. PD)  then the PPS method turns into the  HSS \cite{bai242003} (resp. PSS \cite{bai262005}) method. If $ {\mathcal P}_1$ is PD and  normal,  i.e., $ {\mathcal P}_1{\mathcal P}_1^*={\mathcal P}_1^*{\mathcal P}_1$, and $ {\mathcal P}_2$ is a skew-Hermitian, then the PPS method reduces to the NSS method \cite{bai142007}. Let $\mathcal{A}=\mathcal H_1+\mathcal H_2+ \mathcal S$, where $\mathcal H_1$ and $\mathcal H_2$ are HPSD and $ \mathcal S$ is skew-Hermitian.  Setting ${\mathcal P}_1=\mathcal H_1$ and ${\mathcal P}_2=\mathcal H_2+ \mathcal S$, we obtain the generalized HSS (GHSS) method \cite{benzi312009}. Therefore, the PPS method is a generalization of the HSS, GHSS, PHSS, NSS, PSS, EPSS and SS methods.

Eliminating ${u}^{k+\frac{1}{2}}$ from \eqref{seq} yields
\begin{align}
	\label{PPS}
	{u}^{k+1}=\Gamma_{\PPS} {u}^k+{c},
\end{align}
where \vspace{-1mm}
\begin{align}
	\label{gama}
	\Gamma_{\PPS}= ( \Sigma + {\mathcal P}_2  )^{-1}  ( \Sigma - {\mathcal P}_1   ) ( \Sigma +  {\mathcal P}_1   )^{-1}  ( \Sigma  -  {\mathcal P}_2 ),
\end{align}
and ${c}=2( \Sigma +  {\mathcal P}_2  )^{-1} \Sigma   ( \Sigma +  {\mathcal P}_1    )^{-1} {b}$.
It is known that the PPS  iteration method is convergent for any initial guess $u^ 0 $,  if and only if  $\rho(\Gamma_{\PPS})<1$ \cite{saad822003}. 

Since the matrix  $\Sigma$ is HPD,  the  PPS iteration matrix can be rewritten as 
\begin{align}\nonumber \label{gtil}
	\Gamma_{\PPS} &=\Sigma^{-\frac{1}{2}}	 (\mathcal I+\tilde{\mathcal P}_2)^{-1}(\mathcal I- \tilde{\mathcal P}_1) (\mathcal I+\tilde{\mathcal P}_1)^{-1}   (\mathcal  I-\tilde{\mathcal P}_2)\Sigma^{\frac{1}{2}}\nonumber\\& =\Sigma^{-\frac{1}{2}} (\mathcal I+\tilde{\mathcal P}_2)^{-1}(\mathcal I+\tilde{\mathcal P}_1)^{-1}  (\mathcal I- \tilde{\mathcal P}_1)  (\mathcal  I-\tilde{\mathcal P}_2)\Sigma^{\frac{1}{2}}\nonumber\\&=
	( \Sigma + {\mathcal P}_2  )^{-1}\Sigma  ( \Sigma +  {\mathcal P}_1  )^{-1}  ( \Sigma - {\mathcal P}_1  )\Sigma ^{-1} ( \Sigma -  {\mathcal P}_2  ),
\end{align}
where   \vspace{-2mm}
\begin{align}
	\label{ptild}
	\tilde{\mathcal P}_i:= \Sigma^{-\frac{1}{2}}  {\mathcal P}_i \Sigma^{-\frac{1}{2}}, ~~\text{for} ~~i=1,2.  
\end{align}
If we define
\begin{align}
	\label{MN}
	{ \mathcal{M}_{\PPS}}:= \frac{1}{2}  ( \Sigma + {\mathcal P}_1  )\Sigma ^{-1} ( \Sigma +  {\mathcal P}_2  )\quad \textrm{and} \quad
	{\mathcal N}_{\PPS}:=\frac{1}{2}  (\Sigma - {\mathcal P}_1  )\Sigma ^{-1} ( \Sigma - {\mathcal P}_2  ),
\end{align}
then  $\mathcal{A}=\mathcal{M}_{\PPS}-\mathcal{N}_{\PPS}$,   the matrix $\mathcal{M}_{\PPS}$ is nonsingular  and  
$ \Gamma_{\PPS}=\mathcal{M}_{\PPS}^{-1} {\mathcal{N}_{\PPS}}
$.
Therefore, if we apply a Krylov subspace method such as GMRES  \cite{saad71986} to approximate the solution of the system \eqref{system}, then the matrix $\mathcal{M}_{\PPS}$ can be considered as a preconditioner to the system. Since the  prefactor  $\frac{1}{2}$ in the preconditioner $\mathcal{M}_{\PPS}$  has no effect on the preconditioned system,  we can take the matrix
\begin{align}
	\label{PPSpre}
	{\mathcal P}_{\PPS}  := ( \Sigma+ {\mathcal P}_1  )\Sigma^{-1} ( \Sigma + {\mathcal P}_2  ),
\end{align}
as the PPS  preconditioner for solving \eqref{system}.

We prove  that if   $\mathcal A$ is PD  and one of the following conditions holds, then the PPS method with any HPD matrix $\Sigma$ is convergent:
\begin{enumerate}[(i)]
	\item    If   one of the  matrices ${\mathcal P}_1$ and $ {\mathcal P}_2 $  is HPSD.
	\item   If  one of the  matrices ${\mathcal P}_1$ and $ {\mathcal P}_2 $   is skew-Hermitian.
	\item  If   one of the matrices  ${\mathcal P}_1$ and $ {\mathcal P}_2 $  is  PD.   
	\item    If   $ \nul({\mathcal P}_1)\cup \nul({\mathcal P}_2+ {\mathcal P}_2^*)=\mathbb{C} ^n$ or  $ \nul({\mathcal P}_2)\cup \nul({\mathcal P}_1+ {\mathcal P}_1^*)=\mathbb{C} ^n$.
	\item  
	If the matrix 
	$ {\mathcal A}+     {\mathcal P}_2^*\Sigma^{-1}{\mathcal P}_1^*\Sigma^{-1}{\mathcal A}$
	is PD.
	\item  If  
	${\mathcal P}_1\Sigma^{-1}{\mathcal P}_2={\mathcal P}_2 \Sigma^{-1}{\mathcal P}_1$ (e.g.  if
	${\mathcal P}_1 {\mathcal P}_2={\mathcal P}_2{\mathcal P}_1$ and $ \Sigma=\alpha \mathcal I $). 
\end{enumerate}

\section{Convergence   of the PPS method}\label{sec:conver}  
Let $ {\mathcal P}\in \mathbb{C}^{n\times n} $.
If $\mathcal{I}+ {\mathcal P}$ is invertible,
then we define  the function $f$ as
\begin{align}\label{f}
	\begin{aligned}
		f( {\mathcal P}):=\norm{( \mathcal I-  {\mathcal P})( \mathcal{I}+ {\mathcal P})^{-1}}.
	\end{aligned}	
\end{align}
Then we have
\begin{align*}
	f( {\mathcal P})& =\max_{x\neq 0}\frac{\norm{( \mathcal{I}-{\mathcal P}) x}}{\norm{(\mathcal{I}+ {\mathcal P})x}} =\max_{x\neq 0} \sqrt{\frac{ {x^*( \mathcal{I}+{\mathcal P}^*{\mathcal P}) x -x^*({\mathcal P}+{\mathcal P}^*) x}}{  x^* (\mathcal{I}+{\mathcal P}^*{\mathcal P})x +x^*({\mathcal P}+{\mathcal P}^*)x}}.
\end{align*}
Therefore, we have  the following remark. 
\begin{remark}
	\label{normp}
	Suppose that $-1\notin \sigma( {\mathcal P})$ and $f( {\mathcal P})$ is defined as in \eqref{f}.
	\begin{enumerate}[~~ 1.]
		\item
		If $ {\mathcal P}$ is skew-Hermitian, then for all $ x\neq 0$ we have $ x^*({\mathcal P}+{\mathcal P}^*)x=0 $  and so $f({\mathcal P})= 1$.
		\item
		If \( \mathcal{P} \) is a normal matrix, then there exist matrices \( U \) and \( D \) such that
		$ \mathcal{P} = U^* \mathcal D U $,
		where \( U \) is a unitary matrix and \( \mathcal D \) is a diagonal matrix whose diagonal entries are the eigenvalues of \( \mathcal{P} \).
		Therefore,
		\begin{align*}
			f( {\mathcal P})&  =  \norm{( \mathcal I-  {U^* \mathcal D U})( \mathcal{I}+ {U^* \mathcal D U})^{-1}}
			\\&=  \norm{( \mathcal I-  {  \mathcal D  })( \mathcal{I}+ {  \mathcal D  })^{-1}} 
			\\&={  \max_{\lambda\in \sigma( {\mathcal P})}}{\left|\frac{1-\lambda}{1+\lambda}\right|}.
		\end{align*}
		\item
		The matrix $ {\mathcal P}$ is  PSD if and only if  
		for all $ x\neq 0$, $ x^*({\mathcal P}+{\mathcal P}^*)x\geq 0 $ if and only if  
		for all $ x\neq 0$, 
	\begin{equation*}
		\frac{ {x^*( \mathcal{I}+{\mathcal P}^*{\mathcal P})x-x^*({\mathcal P}+{\mathcal P}^*) x}}{  x^* (\mathcal{I}+{\mathcal P}^*{\mathcal P})x +x^*({\mathcal P}+{\mathcal P}^*)x}\leq 1
	\end{equation*}
		if and only if  
		$	f({\mathcal P})\leq 1$.
		\item
		The matrix $ {\mathcal P}$ is  PD if and only if  
		for all $ x\neq 0$, $ x^*({\mathcal P}+{\mathcal P}^*)x> 0 $  if and only if  			for all $ x\neq 0$, 
	\begin{equation*}
		\frac{ {x^*( \mathcal{I}+{\mathcal P}^*{\mathcal P})x-x^*({\mathcal P}+{\mathcal P}^*) x}}{  x^* (\mathcal{I}+{\mathcal P}^*{\mathcal P})x +x^*({\mathcal P}+{\mathcal P}^*)x}< 1
	\end{equation*}
		if and only if $	f({\mathcal P})< 1$. 
	\end{enumerate}
\end{remark}

Let  $ 	\tilde{\Gamma}:=(\mathcal I+\tilde{\mathcal P}_2)^{-1}(\mathcal I+\tilde{\mathcal P}_1 )^{-1}(\mathcal I- \tilde{\mathcal P}_1 ) 
(\mathcal  I-\tilde{\mathcal P}_2) $.  By \eqref{gtil}, we have 
$\sigma( \Gamma_{\PPS}) =\sigma(\tilde{\Gamma})$ and  moreover  $  x \in  \ev(\Gamma_{\PPS})  $
if and only if $\Sigma^{\frac{1}{2}} x  \in \ev(\tilde{\Gamma})$.

To facilitate the proof of the theorems, we will work with the sets $\ev(\tilde{\Gamma}) $ and $\sigma(\tilde{\Gamma})$. In the next theorem, we present an upper  bound for $\rho(\Gamma_{\PPS})$.
\begin{theorem}
	\label{boundPPS}
	Suppose that the matrix $\mathcal{A}$   has a splitting  of the form \eqref{orgsplit}.
	Then,		
	for every $\lambda \in \sigma(\Gamma_{\PPS})$, there exists a vector $ x_\lambda \in \ev(\tilde{\Gamma}) $, such that   
	\begin{align}\label{landa}
		|\lambda| \leqslant \min\left\{ m_1\left(x  _\lambda\right),m_2\left( y_\lambda\right)\right\},
	\end{align}
	where $y_\lambda:=(\mathcal I+\tilde{\mathcal P}_1)^{-1}   ( \mathcal I-\tilde{\mathcal P}_2)x_\lambda $ and
	\begin{align}\label{mi}
		m_{3-i}(x):=f(\tilde{\mathcal P}_{3-i})\sqrt{\frac{x^*  (\mathcal I+ \tilde{\mathcal P}_i^* \tilde{\mathcal P}_i) x-x^*(\tilde{\mathcal P}_i+\tilde{\mathcal P}_i^*) x}{x^*  (\mathcal I+ \tilde{\mathcal P}_i^* \tilde{\mathcal P}_i) x+x^*(\tilde{\mathcal P}_i+\tilde{\mathcal P}_i^*) x}}\leq 1, ~~\text{for}~~i=1,2,
	\end{align}
	in which $f$ is defined in \eqref{f}. Hence
	\begin{align}\label{rho}
		\rho(\Gamma_{\PPS}) \leqslant \max_{x_\lambda\in \ev(\tilde{\Gamma})  }\min\left \{ m_1\left(x_\lambda  \right),m_2\left((\mathcal I+\tilde{\mathcal P}_1)^{-1}   (\mathcal I-\tilde{\mathcal P}_2) x_\lambda\right)\right \}\leqslant 1.
	\end{align}
\end{theorem}
\begin{proof}		
	Let $\lambda \in \sigma(\Gamma_{\PPS})$.    
	Therefore
	$\lambda \in \sigma (\tilde{\Gamma})$. Thus,
	there exists a vector $x_\lambda\in \ev(\tilde{\Gamma}) $  such that
	$\tilde{\Gamma}x_\lambda=\lambda   x_\lambda$. 
	So, 
\begin{equation*} 
	(\mathcal I+\tilde{\mathcal P}_1)^{-1}(\mathcal I- \tilde{\mathcal P}_1)  (\mathcal I-\tilde{\mathcal P}_2)x_\lambda=\lambda  (\mathcal I+\tilde{\mathcal P}_2) x_\lambda.
\end{equation*}		 
	Hence
	\begin{align}\label{unitary}
\nonumber		\norm{\lambda  (\mathcal I+\tilde{\mathcal P}_2) x_\lambda}&=\norm{(\mathcal I+\tilde{\mathcal P}_1)^{-1}(\mathcal  I- \tilde{\mathcal P}_1)  (\mathcal I-\tilde{\mathcal P}_2)x_\lambda}\\&\leqslant \norm{(\mathcal I+\tilde{\mathcal P}_1)^{-1}(\mathcal I- \tilde{\mathcal P}_1)} \norm{(\mathcal I-\tilde{\mathcal P}_2)x_\lambda}\nonumber\\&
		=f(\tilde{\mathcal P}_1) \norm{(\mathcal I-\tilde{\mathcal P}_2)x_\lambda},
	\end{align}
	which is equivalent to
\begin{equation*} 
 |\lambda | \leq   f(\tilde{\mathcal P}_1) \frac{\norm{(\mathcal I-\tilde{\mathcal P}_2)x_\lambda}}{\norm{(\mathcal I+\tilde{\mathcal P}_2)x_\lambda} } ,	
\end{equation*}
	which results in 
	$|\lambda| \leqslant m_1\left(x_\lambda\right)$.  
	
	Also, from 	$\tilde{\Gamma}x_\lambda=\lambda  x_\lambda$, we have $(\mathcal I-\tilde{\mathcal P}_2) \tilde{\Gamma} x_\lambda=\lambda (\mathcal I-\tilde{\mathcal P}_2)x_\lambda $ and so
	\begin{align*}
		\norm{\lambda (\mathcal I-\tilde{\mathcal P}_2)x_\lambda}
		&=  \norm{(\mathcal I-\tilde{\mathcal P}_2)(\mathcal  I+\tilde{\mathcal P}_2)^{-1}(\mathcal I- \tilde{\mathcal P}_1) (\mathcal I+\tilde{\mathcal P}_1)^{-1}   (\mathcal I-\tilde{\mathcal P}_2)x_\lambda}
		\\&\leq f(\tilde{\mathcal P}_2)\norm{(\mathcal I- \tilde{\mathcal P}_1) (\mathcal I+\tilde{\mathcal P}_1)^{-1}   (\mathcal I-\tilde{\mathcal P}_2)x_\lambda}.
	\end{align*}
	Let $y_\lambda:=(\mathcal I+\tilde{\mathcal P}_1)^{-1}   (\mathcal I-\tilde{\mathcal P}_2)x_\lambda$. 
	Then, $\norm{\lambda (\mathcal I+\tilde{\mathcal P}_1)y_\lambda} \leq f(\tilde{\mathcal P}_2)\norm{(\mathcal I- \tilde{\mathcal P}_1) y_\lambda}$
	which is equivalent to
\begin{equation*} 
	 |\lambda | \leq   f(\tilde{\mathcal P}_2) \frac{\norm{(\mathcal I-\tilde{\mathcal P}_1)y_\lambda}}{\norm{(\mathcal I+\tilde{\mathcal P}_1)y_\lambda} } ,
\end{equation*}
	and so $|\lambda| \leqslant m_2(y_\lambda)$.  
	
	It is  clear that matrices  $\tilde{\mathcal P}_1$ and $\tilde{\mathcal P}_2$ are  PSD. By Remark \ref{normp},     $m_1(x_\lambda)\leq f(\tilde{\mathcal P}_1)f(\tilde{\mathcal P}_2)  \leq 1$. Similarly,   $m_2(y_\lambda)\leq 1$. Hence, the proof is completed.
\end{proof}	

The condition that at least one of the matrices ${\mathcal P}_1$ and ${\mathcal P}_2$ to be  PD is not necessary for $\rho(\Gamma_{\PPS})<1$, but  under this condition, the convergence of the  PPS method will be guaranteed, because in this case, by Remark \ref{normp}, we have
\begin{align}
	\label{mf}
	m_i\left(x \right)\leq   f(\tilde{\mathcal P}_1)f(\tilde{\mathcal P}_2 )<1, ~~~~\text{for all}~~x\in \mathbb{C}^n ~ \text{and}~~  i=1,2.
\end{align} 
\begin{lemma} [\cite{horn2012}]
	If     ${\mathcal P}_2$ is skew-Hermitian, then ${(\mathcal  I- \tilde{\mathcal P}_2)(\mathcal I+\tilde{\mathcal P}_2)^{-1}}$ is a unitary matrix.
\end{lemma}
If $ {\mathcal P}_1$ 	is  PD and ${\mathcal P}_2=0$  (or ${\mathcal P}_2$ is skew-Hermitian), then    by  \eqref{mf}   and Theorem \ref{boundPPS}, the PPS method  is  convergent  and   turns into the SS (PSS) method.  	It  can be  seen that  in this case  
\begin{align*} 
	\rho (\Gamma_{\PPS}) &
	=\nonumber
	\rho( \tilde{\Gamma} )=   \rho\big(  (\mathcal I-\tilde{\mathcal P}_2)(\mathcal  I+\tilde{\mathcal P}_2) ^{-1}(\mathcal I+\tilde{\mathcal P}_1)^{-1}(\mathcal  I- \tilde{\mathcal P}_1)\big)\\&=
	\max_{x \in \ev(\tilde{\Gamma}) }\frac{\norm{\big(  (\mathcal I-\tilde{\mathcal P}_2)(\mathcal  I+\tilde{\mathcal P}_2) ^{-1}(\mathcal I+\tilde{\mathcal P}_1)^{-1}(\mathcal  I- \tilde{\mathcal P}_1)\big)x}}{\norm{ x}}
	\\&=
	\max_{x \in \ev(\tilde{\Gamma}) }\frac{\norm{\big(  (\mathcal I+\tilde{\mathcal P}_1)^{-1}(\mathcal  I- \tilde{\mathcal P}_1)\big)x}}{\norm{ x}}
	\end{align*}
which provides a value   for $ 	\rho (\Gamma_{\PPS})  $.
Moreover, if $\tilde{\mathcal P}_1$ is normal,  then    by Remark \ref{normp},  we have 
\begin{align*} 
	\rho(\Gamma_{\PPS})&\leq  
	\max_{x  \neq 0}\frac{\norm{\big(  (\mathcal I+\tilde{\mathcal P}_1)^{-1}(\mathcal  I- \tilde{\mathcal P}_1)\big)x}}{\norm{ x}}= f(\tilde{\mathcal P}_1)\\&=\max_{\mu\in \sigma(\tilde{\mathcal P}_1)} {\frac{\left| 1-\mu\right | }{\left|1+\mu\right|}} = \max_{\mu\in \sigma(\tilde{\mathcal P}_1)} \sqrt{{\frac{(1+\mu\bar{\mu})-(\mu+\bar{\mu})}{(1+\mu\bar{\mu})+(\mu+\bar{\mu})}}}<1,
\end{align*}
In this case, the PPS method reduces to HSS or NSS method  (see \cite[Theorem  2.2]{bai242003} and \cite[Theorem 2.2]{bai142007}).

In Theorem \ref{boundPPS},  we gave  an upper bound for the spectral  radius of $\Gamma_{\PPS}$. 
In the  following, we obtain the spectral  radius of $\Gamma_{\PPS}$ and a necessary and sufficient condition for the convergence of the PPS method.
\begin{theorem}
	\label{convPPS}
	Suppose that the matrix $\mathcal{A}$    has a splitting  of the form \eqref{orgsplit}  and  
	\begin{align}\label{r_i}
		\begin{cases}
			r_1(x):=x^*\left(\mathcal I+\tilde{\mathcal A}^*\tilde{\mathcal A}+\tilde{\mathcal P}_1\tilde{\mathcal P}_2+\tilde{\mathcal P}_2^*\tilde{\mathcal P}_1^*+\tilde{\mathcal P}_2^*\tilde{\mathcal P}_1^* \tilde{\mathcal P}_1\tilde{\mathcal P}_2\right)x,\\
			r_2(x):=x^*\left(\tilde{\mathcal A}+\tilde{\mathcal A}^*+\tilde{\mathcal A}^*\tilde{\mathcal P}_1\tilde{\mathcal P}_2+ \tilde{\mathcal P}_2^*\tilde{\mathcal P}_1^*\tilde{\mathcal A}\right)x,
		\end{cases}
	\end{align} 
	where $\tilde{\mathcal A}=\Sigma^{-\frac{1}{2}} {\mathcal A} \Sigma^{-\frac{1}{2}}$. Then,  for all $x_\lambda \in \ev(\tilde{\Gamma}) $, we have $   r_1(x_\lambda)\geq r_2(x_\lambda)\geq 0$,  $r_1(x_\lambda)>0$
	and 
	\begin{equation*} 
		\rho(\Gamma_{\PPS})=\max_{x_\lambda \in \ev(\tilde{\Gamma}) }
		\sqrt{\frac{r_1(x_\lambda)-r_2(x_\lambda)}{ r_1(x_\lambda)+r_2(x_\lambda)}}.
	\end{equation*}
	So,  $\rho(\Gamma_{\PPS})<1$ if and only if   for all $x_\lambda \in \ev(\tilde{\Gamma}) $, we have   $   {r_2}(x_\lambda)\neq 0 $.
\end{theorem}
\begin{proof}
	Let $\lambda \in \sigma(\Gamma_{\PPS})$. By Theorem \ref{boundPPS},   there exists a vector $0 \neq x_\lambda\in \mathbb{C}^{n}$, such that
	\begin{align} 
		(\mathcal I-\tilde{\mathcal P}_1)(\mathcal I-\tilde{\mathcal P}_2) x_\lambda=\lambda (\mathcal I+\tilde{\mathcal P}_1)(\mathcal I+\tilde{\mathcal P}_2)  x_\lambda.
	\end{align} 
	Therefore, $\norm{(\mathcal  I-\tilde{\mathcal A}+\tilde{\mathcal P}_1\tilde{\mathcal P}_2) x_\lambda }=\norm{\lambda (\mathcal I+\tilde{\mathcal A}+\tilde{\mathcal P}_1\tilde{\mathcal P}_2)  x_\lambda}$.
	Hence, 
	\begin{align*}
		|\lambda|^2&=\frac{\norm{(\mathcal  I-\tilde{\mathcal A}+\tilde{\mathcal P}_1\tilde{\mathcal P}_2) x_\lambda }^2}{\norm{(\mathcal I+\tilde{\mathcal A}+\tilde{\mathcal P}_1\tilde{\mathcal P}_2)  x_\lambda}^2} = {\frac{r_1(x_\lambda)-r_2(x_\lambda)}{ r_1(x_\lambda)+r_2(x_\lambda) }},
	\end{align*}
	where    $	r_i(x)$, for $i=1,2$,  are defined in  \eqref{r_i}.
	By \eqref{landa} and \eqref{mi}, 
\begin{equation*} 
{\frac{r_1(x_\lambda)-r_2(x_\lambda)}{ r_1(x_\lambda)+r_2(x_\lambda) }}=|\lambda|^2\leq (m_1(x_\lambda) )^2\leq 1	
\end{equation*}
	and so $r_2(x_\lambda)\geq 0 .$
	Therefore,  from  
\begin{equation*} 
	r_1(x_\lambda)-r_2(x_\lambda)=\norm{(\mathcal I-\tilde{\mathcal P}_1)(\mathcal I-\tilde{\mathcal P}_2) x_\lambda}^2\geq 0,
\end{equation*}
	we deduce that $r_1(x_\lambda)\geq r_2(x_\lambda) $.
	Also, since $(\mathcal I+\tilde{\mathcal P}_1)(\mathcal I+\tilde{\mathcal P}_2) $ is nonsingular, we get $$r_1(x_\lambda)+r_2(x_\lambda)=\norm{(\mathcal I+\tilde{\mathcal P}_1)(\mathcal I+\tilde{\mathcal P}_2) x_\lambda}^2> 0,$$ and so $r_1(x_\lambda)>0$.  Hence 
	$$ 	\rho(\Gamma_{\PPS})=\max_{x_\lambda \in \ev(\tilde{\Gamma}) }
	\sqrt{\frac{r_1(x_\lambda)-r_2(x_\lambda)}{ r_1(x_\lambda)+r_2(x_\lambda)}},~~   r_1(x_\lambda)\geq r_2(x_\lambda)\geq 0, ~~r_1(x_\lambda)>0, $$ 
	and 
	$\rho(\Gamma_{\PPS})<1$ if and only if  $ r_2(x_\lambda)\neq 0$, for all $x_\lambda \in \ev(\tilde{\Gamma}) $.
\end{proof}

 Using Theorem \ref{convPPS}, we have following corollary. 
\begin{corollary}\label{convK}
	Suppose that the matrix $\mathcal{A}$   has a splitting  of the form \eqref{orgsplit} and $K=\tilde{\mathcal A}+ \tilde{\mathcal P}_2^*\tilde{\mathcal P}_1^*\tilde{\mathcal A}$. If $K$ 
	is	PD,  then the PPS method is convergent.
\end{corollary}
\begin{proof}
	Suppose that the matrix $K$   is PD.  Therefore, $K+K^*$ is HPD and so by  \cite[Corollary 4.3.15]{horn2012}, $
	r_2(x)=x^*(K+K^*)x\geq 
	\lambda_{\min} \left(K+K^*\right) \|x\|_2^2>0,
	$  for all $x \in \mathbb{C}^n$. Thus, by Theorem \ref{convPPS}, $\rho(\Gamma_{\PPS})<1.$
\end{proof}
\begin{corollary}
	Suppose that the matrix $\mathcal{A}$ is  PD and   has a splitting  of the form \eqref{orgsplit}. If    $\tilde{\mathcal P}_1\tilde{\mathcal P}_2=\tilde{\mathcal P}_2\tilde{\mathcal P}_1$,  then the PPS method is convergent.  
\end{corollary}
\begin{proof}
	Since    $\tilde{\mathcal P}_1\tilde{\mathcal P}_2=\tilde{\mathcal P}_2\tilde{\mathcal P}_1$, we have 
	\begin{align*}
		K+K^*&=  \tilde{\mathcal A}+\tilde{\mathcal A}^*+\tilde{\mathcal A}^*\tilde{\mathcal P}_1\tilde{\mathcal P}_2+ \tilde{\mathcal P}_2^*\tilde{\mathcal P}_1^*\tilde{\mathcal A}
		\\&= \tilde{\mathcal A}+\tilde{\mathcal A}^*+(\tilde{\mathcal P}_1+\tilde{\mathcal P}_2)^*\tilde{\mathcal P}_1\tilde{\mathcal P}_2+ \tilde{\mathcal P}_2^*\tilde{\mathcal P}_1^*(\tilde{\mathcal P}_1+\tilde{\mathcal P}_2) \\
		&= \tilde{\mathcal A}+\tilde{\mathcal A}^*+\tilde{\mathcal P}_2^*(\tilde{\mathcal P}_1+\tilde{\mathcal P}_1^*)\tilde{\mathcal P}_2+\tilde{\mathcal P}_1^*\tilde{\mathcal P}_1\tilde{\mathcal P}_2+ \tilde{\mathcal P}_2^*\tilde{\mathcal P}_1^*\tilde{\mathcal P}_1 \\
		&= \tilde{\mathcal A}+\tilde{\mathcal A}^*+\tilde{\mathcal P}_2^*(\tilde{\mathcal P}_1+\tilde{\mathcal P}_1^*)\tilde{\mathcal P}_2+\tilde{\mathcal P}_1^*(\tilde{\mathcal P}_2+\tilde{\mathcal P}_2^*)\tilde{\mathcal P}_1.
	\end{align*}
	Since   ${\mathcal A}$ is  PD, we conclude that $ \tilde{\mathcal A}+\tilde{\mathcal A}^ *$ is HPD. Therefore $K+K^*$ is HPD and so
	$K $ is PD. Hence, by Corollary  \ref{convK}, the PPS  method is convergent. 
\end{proof}

The next theorem presents other conditions for the convergence of the PPS method. 

\begin{theorem}\label{otherconv}
	Suppose that the matrix $\mathcal{A}$ is PD and  has a splitting  of the form \eqref{orgsplit}.  If one of the following conditions holds, then $\rho(\Gamma_{\PPS})< 1$:
	\begin{enumerate}[(1)]
		\item \label{cond:1}
		If   $  \ev(\Gamma_{\PPS}) \cap    \nul{({\mathcal P}_2+ {\mathcal P}_2^*)}\subseteq    \nul({\mathcal P}_2)$.
		\item  \label{cond:2}
		If   $  \ev(\Gamma_{\PPS})  \subseteq    \nul({\mathcal P}_2)\cup  \nul({\mathcal P}_1+ {\mathcal P}_1^*)$.
	\end{enumerate}
\end{theorem}

\begin{proof}  Let $(\lambda,x)$ be an eigenpair  of the matrix $\Gamma_{\PPS}$, i.e., $\Gamma_{\PPS}x= \lambda  x $.  Thus, $  x\in \ev(\Gamma_{\PPS}) $ and  we get $ x_\lambda=\Sigma ^\frac{1}{2}x \in \ev(\tilde{\Gamma})  $.  Since $\mathcal A$ is PD, $x ^* ({\mathcal A}+{\mathcal A}^* )x>0$ and so    $x ^* ({\mathcal P}_2+{\mathcal P}_2^* )x>0$ or $x ^* ({\mathcal P}_1+{\mathcal P}_1^* )x>0$.
	\begin{enumerate}[(1)]
		\item[(1)]  If $x ^* ({\mathcal P}_2+{\mathcal P}_2^* )x>0$, then $x_\lambda^* (\tilde{\mathcal P}_2+\tilde{\mathcal P}_2^* )x_\lambda>0$ and as a result by Theorem \ref{boundPPS}, we have $ m_1\left(x_\lambda \right)<1$ and so  $|\lambda|<1 $. On the other hand, if $x ^* ({\mathcal P}_2+{\mathcal P}_2^* )x=0$, then $x ^* ({\mathcal P}_1+{\mathcal P}_1^* )x>0$. 
		Moreover, by   \cite[Observation 7.1.6]{horn2012}, we have  $x\in \nul({\mathcal P}_2+{\mathcal P}_2^*)$ and so  $ x \in  \ev(\Gamma_{\PPS}) \cap    \nul{({\mathcal P}_2+ {\mathcal P}_2^*)}$. 
		By the assumption,  $x \in  \nul({\mathcal P}_2)$ and hence  $ {\mathcal P}_2x=0 $.  Therefore, $ \tilde{\mathcal P}_2x_\lambda=0 $.
		 {So  we have $x_\lambda^*\tilde{\mathcal P}_2^*=0$. By the definition  of  $r_2$ in   \eqref{r_i}, we conclude that 
			\begin{align*}
				r_2(x_\lambda)=~&x_\lambda^* ( \tilde{\mathcal A}+\tilde{\mathcal A}^*+\tilde{\mathcal A}^*\tilde{\mathcal P}_1\tilde{\mathcal P}_2+ \tilde{\mathcal P}_2^*\tilde{\mathcal P}_1^*\tilde{\mathcal A} )x_\lambda\\
				=
				~&x_\lambda^* (\tilde{\mathcal P}_1+\tilde{\mathcal P}_1^*+\tilde{\mathcal P}_2+ \tilde{\mathcal P}_2^* )x_\lambda+x_\lambda^* (\tilde{\mathcal A}^*\tilde{\mathcal P}_1)\tilde{\mathcal P}_2 x_\lambda+ x_\lambda^* \tilde{\mathcal P}_2^*(\tilde{\mathcal P}_1^*\tilde{\mathcal A} )x_\lambda\\
				=~&x_\lambda^* (\tilde{\mathcal P}_1+\tilde{\mathcal P}_1^*)x_\lambda\\
				=~&x ^* ({\mathcal P}_1+{\mathcal P}_1^* )x>0.
		\end{align*}} 
		So,   by Theorem \ref{convPPS},  $ \rho(\Gamma_{\PPS})<1 $.
		\item [(2)] 
		Let $  \ev(\Gamma_{\PPS})  \subseteq    \nul({\mathcal P}_2)\cup  \nul({\mathcal P}_1+ {\mathcal P}_1^*)$. 	 Since, $ \mathcal A $ is  PD, we have $  \nul({\mathcal P}_1+ {\mathcal P}_1^*)\cap      \nul{({\mathcal P}_2+ {\mathcal P}_2^*)}= \{0\}$ and so
		\begin{align*}
			\ev(\Gamma_{\PPS})  \cap    \nul{({\mathcal P}_2+ {\mathcal P}_2^*)}& \subseteq ( \nul({\mathcal P}_2)\cup  \nul({\mathcal P}_1+ {\mathcal P}_1^*))  \cap    \nul{({\mathcal P}_2+ {\mathcal P}_2^*)}\\&
			= ( \nul({\mathcal P}_2) \cap    \nul{({\mathcal P}_2+ {\mathcal P}_2^*)})\cup \{0\}  .
		\end{align*}
		Obviously $ \nul({\mathcal P}_2) \subseteq    \nul{({\mathcal P}_2+ {\mathcal P}_2^*)} $.  Therefore
		\begin{align*}
			\ev(\Gamma_{\PPS})  \cap    \nul{({\mathcal P}_2+ {\mathcal P}_2^*)}& \subseteq 
			\nul({\mathcal P}_2)   .
		\end{align*}
		Using case (1),  we have  $\rho(\Gamma_{\PPS})< 1$.  \qedhere	\end{enumerate} 
\end{proof}
Using Theorem \ref{otherconv}, we have the following corollary. 
\begin{corollary} \label{cor8}
	Let $\mathcal{A}$ be a PD matrix with a splitting  of the form \eqref{orgsplit}. 
	If  ${\mathcal P}_2$ is  PD,  or  ${\mathcal P}_2$ is HPSD, or ${\mathcal P}_1$ is skew-Hermitian, or   $ \nul({\mathcal P}_2)\cup \nul({\mathcal P}_1+ {\mathcal P}_1^*)=\mathbb{C} ^n$, then $\rho(\Gamma_{\PPS})<1$.
\end{corollary}
 {\begin{proof}
		If $\mathcal{P}_2$ is PD, then $\nul{(\mathcal{P}_2 + \mathcal{P}_2^*)} = \{0\}$, and so
	\begin{equation*} 
		  \ev(\Gamma_{\PPS}) \cap    \nul{({\mathcal P}_2+ {\mathcal P}_2^*)} = \{0\}\subseteq    \nul({\mathcal P}_2).
	\end{equation*}
		Hence,  the result follows from  part  (1) of  Theorem \ref{otherconv}.
		
		If $\mathcal{P}_2$ is HPSD, then $\nul{(\mathcal{P}_2 + \mathcal{P}_2^*)} = \nul{(2\mathcal{P}_2)} = \nul{(\mathcal{P}_2)}$, and 
	\begin{equation*} 
		 \ev(\Gamma_{\PPS}) \cap    \nul{({\mathcal P}_2+ {\mathcal P}_2^*)} = \ev(\Gamma_{\PPS}) \cap    \nul{({\mathcal P}_2)} \subseteq    \nul({\mathcal P}_2).
	\end{equation*}
	Therefore, convergence of the method is deduced using   part (1) of  Theorem \ref{otherconv}.

		If $\mathcal{P}_1$ is skew-Hermitian, then $\nul{(\mathcal{P}_1 + \mathcal{P}_1^*)} = \nul(0) = \mathbb{C}^n$, and hence
	\begin{equation*} 
		  \ev(\Gamma_{\PPS})  \subseteq \mathbb{C}^n= \nul({\mathcal P}_2)\cup  \nul({\mathcal P}_1+ {\mathcal P}_1^*).
	\end{equation*}
		Thus, convergence of the method is concluded using  part (2) of Theorem \ref{otherconv}. 
		
		If $ \nul({\mathcal P}_2)\cup \nul({\mathcal P}_1+ {\mathcal P}_1^*)=\mathbb{C} ^n$,  then
		convergence of the method is inferred from part (2) of  Theorem \ref{otherconv}. 
		
\end{proof}}

\begin{remark}\label{Rem3.9}
	We have 
	\begin{align*}
		\sigma(\Gamma_{\PPS})&=	\sigma\Big((\mathcal I+\tilde{\mathcal P}_2)^{-1}(\mathcal I+\tilde{\mathcal P}_1)^{-1}   (\mathcal I- \tilde{\mathcal P}_1)  (\mathcal  I-\tilde{\mathcal P}_2)\Big)
		\\&=	\sigma\left((\mathcal  I-\tilde{\mathcal P}_2)(\mathcal I+\tilde{\mathcal P}_2)^{-1}(\mathcal I+\tilde{\mathcal P}_1)^{-1}  (\mathcal I- \tilde{\mathcal P}_1)  \right)
		\\&=	\sigma\left((\mathcal I+\tilde{\mathcal P}_2)^{-1}(\mathcal  I-\tilde{\mathcal P}_2)(\mathcal I- \tilde{\mathcal P}_1) (\mathcal I+\tilde{\mathcal P}_1)^{-1}   \right)
		\\&=	\sigma\left((\mathcal I+\tilde{\mathcal P}_1)^{-1} (\mathcal I+\tilde{\mathcal P}_2)^{-1}(\mathcal  I-\tilde{\mathcal P}_2)(\mathcal I- \tilde{\mathcal P}_1)   \right)
	\end{align*} 
	Therefore, by changing 
	the roles of the matrices ${\mathcal P}_1$ and  ${\mathcal P}_2$, i.e.,   considering the splitting  $\mathcal{A}= {\mathcal P}_2+ {\mathcal P}_1$,   
	if the matrix $\mathcal{A}$ is PD and one of the following conditions holds,   then the PPS method with any HPD matrix $\Sigma$ is convergent: 
	\begin{enumerate}[(i)]
		\item    If   one of the  matrices ${\mathcal P}_1$ and $ {\mathcal P}_2$  is HPSD.
		\item   If  one of the  matrices ${\mathcal P}_1$ and $ {\mathcal P}_2 $  is skew-Hermitian.
		\item  If   one of the matrices  ${\mathcal P}_1$ and $ {\mathcal P}_2 $ is  PD.   
		\item  \label{nulp1p2} If   $ \nul({\mathcal P}_1)\cup \nul({\mathcal P}_2+ {\mathcal P}_2^*)=\mathbb{C} ^n$ or  $ \nul({\mathcal P}_2)\cup \nul({\mathcal P}_1+ {\mathcal P}_1^*)=\mathbb{C} ^n$.
	\end{enumerate}
\end{remark}

Now, some comments and suggestions are in order.  In each of the following cases  the PPS method is  convergent:   
\begin{enumerate} [(1)]
	\item $\mathcal  A=\mathcal H+{\mathcal P}$ is  PD, where $\mathcal H$ is HPSD and ${\mathcal P}$ is PSD.
	
	\item  $\mathcal  A={\mathcal P}_1+{\mathcal P}_2$ and  ${\mathcal P}_1$  is PD  and ${\mathcal P}_2$ is PSD.

	\item  $\mathcal  A={\mathcal P}+\mathcal S$ is PD, where ${\mathcal P}$ is PSD and  $\mathcal S$ is skew-Hermitian. 
	It is worth noting that in the special case that  $\Sigma=\alpha \mathcal I$, the HSS, NSS, PSS,  EPSS, and SS iterative methods are obtained. 
	
\end{enumerate}

Consider the saddle point matrix
\begin{equation}\label{Eqq0}
	\mathcal  A=\begin{bmatrix}
		A&B^*\\-B&C
	\end{bmatrix},
\end{equation}
where $A\in \mathbb{C}^{n\times n}$,  $B\in \mathbb{C}^{m\times n}$  and $C\in \mathbb{C}^{m\times m}$. We assume that the matrix    $\mathcal  A$ is PD (in this case, the matrices $A$ and $C$ are PD).  Then, from \eqref{nulp1p2} of  Remark \ref{Rem3.9},  the PPS method is convergent with the splitting 
$\mathcal  A={\mathcal P}_1+{\mathcal P}_2$, in which
\begin{equation}\label{EqqA}
	{\mathcal P}_1=\begin{bmatrix}
		0&0\\0&C
	\end{bmatrix},\quad 
	{\mathcal P}_2=\begin{bmatrix}
		A&B^*\\-B&0
	\end{bmatrix}, 
\end{equation}
or
\begin{equation}\label{EqqB}
	{\mathcal P}_1=\begin{bmatrix}
		A&0\\0&0
	\end{bmatrix}, \quad
	{\mathcal P}_2=\begin{bmatrix}
		0&B^*\\-B&C
	\end{bmatrix}.
\end{equation}
However, when the matrix  $\mathcal A$ is PSD,   convergence of the PPS method is not guaranteed.
For instance, consider the nonsingular matrix
$
\mathcal{A}=\left[\begin{smallmatrix}
	0&-1\\
	1&0
\end{smallmatrix}\right],
$
with splitting $\mathcal{A}= {\mathcal P}_1+ {\mathcal P}_2$, in which  $ {\mathcal P}_1=\mathcal{A}$ and $ {\mathcal P}_2=0$. Letting $\Sigma=\mathcal I$, we have $\Gamma _{\PPS}=\mathcal{A}$. The eigenvalues of the iteration matrix $\Gamma _{\PPS}$ are  $\lambda=\pm  i$.
So, the spectral radius of $\Gamma _{\PPS}$ is equal to $1$ and  the PPS method is not convergent. 
Therefore, if the matrix   $\mathcal{A}$ is nonsingular and the convergence conditions are not met, then the PPS method may diverge.
Nevertheless,    we can  generate a new convergent  iterative method for nonsingular    PSD  system  \eqref{system}  as (see \cite[page 27]{benzi262004}, \cite{Song1061999})
\begin{align}
	\label{iter2}
	\hat u^{k+1}=\left((1-\beta )I+\beta \Gamma_{\PPS}\right)\hat u^k +\beta \mathcal M_{\PPS} ^{-1} b,~~~ \beta \in (0,1),~~k=1,2,\ldots
\end{align}
where $\mathcal M_{\PPS}$  is defined in \eqref{MN} and $ \hat u^{0} $ is initial guess.
Therefore, we can use the preconditioner ${\mathcal P}_\beta=\frac{1}{\beta} \mathcal M_{\PPS}=\frac{1}{2\beta} {\mathcal P}_{\PPS}$ for the  equation  \eqref{system}, for all $\beta \in (0,1)$, where  ${\mathcal P}_{\PPS}$ is presented in  \eqref{PPSpre}.
It is advisable to use a value of $ \beta $ close to 1.
Since $\beta $ is close to  1 and  the  prefactor  $\frac{1}{2\beta}$   has no effect on the preconditioned system, 
we can continue using the preconditioner ${\mathcal P}_{\PPS}$ even if none of the provided convergence conditions are met.

\section{Choosing the  matrix $\Sigma$}\label{sec:sigma}\label{choice sigma}

If $ {\mathcal P}_1$ is HPD   and $\Sigma= {\mathcal P}_1$, then  $\rho(\Gamma _{\PPS})=\rho(0)=0$. 
However, in this scenario, solving the first half-step of the system in \eqref{seq} is equivalent to solving $\mathcal{A}{u}={b}$. Thus, selecting a matrix $\Sigma$  close to ${\mathcal P}_1$ can lead to faster convergence with the PPS method.
This is because when $\Sigma$ is close to ${\mathcal P}_1$, the matrix $\tilde{\mathcal P}_1$ is also close to the identity matrix,  implying that   $m_1$ in Theorem \ref{boundPPS}  approaches zero. Consequently,   $\rho(\Gamma_{\PPS}) \approx 0$, indicating rapid convergence of the PPS method. However, this may increase the computation time.

Since the matrix ${\mathcal P}_1$ is   PSD, we see that the matrix $\mathcal{H}_1=\frac{1}{2} ({\mathcal P}_1+{\mathcal P}_1^*)$ is   HPSD. If $\mathcal{H}_1$ is nonsingular, then $\Sigma$ can be chosen as an approximation of $\mathcal{H}_1$. Otherwise, $\Sigma$ can be selected as an approximation of the matrix $\epsilon I+\mathcal{H}_1$, where $\epsilon$ is a small positive number.
The same relationship holds for ${\mathcal H}_2=\frac{1}{2} ({\mathcal P}_2+{\mathcal P}_2^*)$. 
So, we take 	$\Sigma$  as an approximation of the matrix 
\begin{align}\label{RemS}
	\frac{1}{2}({\mathcal H}_1+{\mathcal H}_2) =\frac{1}{4} ({\mathcal A}+{\mathcal A}^*) .  
\end{align}

In the sequel, we assume that $\Sigma=\alpha Q$, where $\alpha>0$ and $ Q $ is a suitable HPD matrix. We choose the parameter $\alpha$ such a way that $ {\mathcal N}_{\PPS}\approx 0 $. In this case, we see that $ \mathcal{A}\approx\mathcal{M}_ {\PPS}$, and a high speed of convergence is expected. To do so, we set
\[
\frac{1}{2}  (\Sigma    - {\mathcal P}_1  )\Sigma ^{-1} ( \Sigma - {\mathcal P}_2  ) =\frac{1}{2}  (\alpha Q  - {\mathcal P}_1  )(\alpha Q ) ^{-1} ( \alpha Q - {\mathcal P}_2  ) \approx 0,   
\]
which results in
\begin{align}\label{optalpha}
	\alpha^2Q+ {\mathcal P}_1Q^{-1}{\mathcal P}_2\approx\alpha {\mathcal A} .
\end{align}
Now, we define the function
$$\phi(\alpha)= \alpha^2 \norm{Q}_F  + \norm{{\mathcal P}_1Q^{-1}{\mathcal P}_2}_F-\alpha    \norm{ {\mathcal A} }_F    ,$$
where  $\norm{\cdot} _F $   is the  Frobenius norm. This function has a minimum at 
\begin{align}\label{alpha2star}
	\alpha_{*}=\frac{\norm{ {\mathcal A}}_F}{2\norm{Q}_F}. 
\end{align}
For this value of $\alpha$, we have 
\begin{equation*} 
	\Sigma=\frac{\norm{ {\mathcal A}}_F}{2} \frac{Q}{\norm{Q}_F},
\end{equation*}
which demonstrates that if we replace the matrix  $Q$ by $kQ$ where $k$ is a positive number, then the matrix $\Sigma$ remains unchanged. Therefore, referring to \eqref{RemS}, we can set $\Sigma=\alpha_{*} Q$, where $\alpha_{*}$  is   defined  by \eqref{alpha2star} and $Q$ is  an approximation of the matrix $ {\mathcal H}=\frac{1}{2} ({\mathcal A}+{\mathcal A}^*) $.

\section{Choosing matrices ${\mathcal P}_1$ and ${\mathcal P}_2$}\label{sec:p1p2}  
For a   PSD matrix $\mathcal{A}$, 
some choices of  ${\mathcal P}_1$ and ${\mathcal P}_2$ are as follows:
\begin{enumerate}[(i)]
	\item  Suppose  that the matrix $\mathcal{A}$	has a  splitting of the form $\mathcal{A}=\mathcal{L}+\mathcal{D}+\mathcal{U}$, where $\mathcal{L}$ and $\mathcal{U}^*$ are  block lower triangular, and $\mathcal{D}$ is block diagonal matrix. Moreover, we assume that the matrix $D$ is PD and $\mathcal{L}=-\mathcal{U}^*$.  Let also $\mathcal{D}=\mathcal{D}_1+\mathcal{D}_2$, where  $\mathcal{D}_1$ and $\mathcal{D}_2$ are positive  semidefinite with $\nul(\mathcal{D}_1)\cup  \nul(\mathcal{D}_2+ \mathcal{D}_2^*)=\mathbb{C}^n$. 
	In this  case, a choice for the matrices ${\mathcal P}_1$ and ${\mathcal P}_2$ is
	\begin{align}
		{\mathcal P}_1=\mathcal{D}_1,
		~~	{\mathcal P}_2 =\mathcal{L}+\mathcal{D}_2+\mathcal{U}.
	\end{align}
	Thus,    ${\mathcal P}_1$ and ${\mathcal P}_2$ are  PSD. Since $\nul({\mathcal P}_1)\cup  \nul({\mathcal P}_2+ {\mathcal P}_2^*)=\mathbb{C}^n$,   by Theorem \ref{boundPPS},    $\rho(\Gamma_{\PPS})< 1$. 
	
	\item  Let  $\mathcal{A}=\hat{\mathcal P}_1+\hat{\mathcal P}_2$ be a PP splitting of the  matrix $\mathcal{A}$, where $\hat{\mathcal P}_1=\hat{\mathcal P}_3+\hat{\mathcal P}_4$ in which $\hat{\mathcal P}_3$ and  $\hat{\mathcal P}_4$ are  PSD. Then, we  can choose ${\mathcal P}_1=\hat{\mathcal P}_3$ and ${\mathcal P}_2=\hat{\mathcal P}_2+\hat{\mathcal P}_4$. 
	For example, let $\mathcal{A}=\mathcal{H}+\mathcal{S}$ be the 
	Hermitian and skew-Hermitian splitting of the matrix 
	$\mathcal{A}$ and  $\mathcal{H}=\mathcal{G}+\mathcal{K}$ (see \cite{benzi312009}), where $\mathcal{K}$ is of a simple
	form (e.g., diagonal). In this case, the matrix $\mathcal{K}$ can be associated with the skew-Hermitian portion $\mathcal{S}$  of $\mathcal{A}$. Hence, we can choose 	${\mathcal P}_1=\mathcal{G}$ and ${\mathcal P}_2=\mathcal{K}+\mathcal{S}$. Assuming that the matrix $\mathcal{A}$ is  PD,
	by applying Theorem \ref{otherconv}, if    $\mathcal{G}$ is HPSD  or $\mathcal{K}$  is skew-Hermitian or  one of the matrices $\mathcal{G}$ and $\mathcal{K}$ is  PD, then for every HPD matrix $\Sigma$, we have $\rho(\Gamma_{\PPS})< 1$ which ensures the convergence of the PPS method. It is important to note that this case represents a generalization of the GHSS method discussed in \cite{benzi312009}. 
	
	\end{enumerate}


\section{Numerical examples}\label{sec:examp}

 {	
	We	consider  the system of linear equations    
	\begin{align}\label{sad}
		\mathcal{A}u=\begin{bmatrix}
			A&B\\C&D
		\end{bmatrix}\begin{bmatrix}
			y_1\\y_2
		\end{bmatrix}=\begin{bmatrix}
			x_1\\x_2
		\end{bmatrix},
	\end{align}
	where  $ \mathcal{A} $ is  PSD. Suppose that  $\mathcal D_\mathcal {A}$, $\mathcal L_\mathcal {A}$ and $\mathcal U_\mathcal {A}$ are the diagonal, strictly  lower triangular and strictly upper triangular matrices of the matrix  $\mathcal{A}$, respectively,  and  
	\begin{align*}
		&\mathcal {BD}_\mathcal {A}=\begin{bmatrix}
			A&0\\0&D
		\end{bmatrix}, 
		~~
		\mathcal {BL}_\mathcal {A}=\begin{bmatrix}
			0&0\\C&0
		\end{bmatrix},
		~~ \mathcal {BU}_\mathcal {A}=\begin{bmatrix}
			0&B\\0&0
		\end{bmatrix}.
	\end{align*}
	We now consider the splitting  $\mathcal{A}={\mathcal P}_1+{\mathcal P}_2$,   with the following choices for the matrix  $ {\mathcal P}_1  $:
	\begin{enumerate}[(1)]
		
		\item     ${\mathcal P}_1=\mathcal D_\mathcal {A}+\mathcal L_\mathcal {A}+\mathcal U_\mathcal {A}^*$.  
		
		In this case, if $\Sigma=\alpha \mathcal   I$,   then the  PPS method coincides with  the triangular and skew-Hermitian splitting (TSS) \cite{bai262005}. If the matrix $\Sigma$ is arbitrary, we call the iteration method as extended TSS (ETSS).
		
		\item  ${\mathcal P}_1=\mathcal {BD}_\mathcal {A}+\mathcal {BL}_\mathcal {A}+\mathcal {BU}_\mathcal {A}^*$.
		
		In this case, if $\Sigma=\alpha \mathcal   I$, then the  PPS method corresponds to the block triangular and skew-Hermitian splitting (BTSS)  \cite{bai262005}. For an arbitrary matrix $\Sigma$, we call the corresponding splitting as extended BTSS (EBTSS).
		Note that, if 	$ C=	 -   B^*$,  then the BTSS method  turns into the PSS \cite{pan1722006} iteration method. 
		
		\item  ${\mathcal P}_1=\frac{1}{2} (\mathcal {A}+\mathcal {A}^*) $.
		
		In this case, if we set $\Sigma=\alpha \mathcal   I$, then the PPS method reduces to the  HSS method \cite{bai242003}.
		If the matrix $\Sigma$ is arbitrary,  we call the corresponding splitting as extended HSS (EHSS).

		\item 	 ${\mathcal P}_1=\mathcal {A}$.
		
		When $\Sigma=\alpha \mathcal   I$, the PPS iteration method concurs with the SS method. For  $\Sigma$ being arbitrary,  we call the corresponding splitting as extended SS (ESS). 
\end{enumerate}}

 {Each of the stationary iteration methods generated using the above splittings to solve the system \eqref{system} induces a preconditioner to the system.}

 {In order to improve the performance of the PPS method, we consider  two special cases of the PP splitting   of the matrix $\mathcal{A} $    as $\mathcal{A}={\mathcal P}_1+{\mathcal P}_2$ with  the HPD  matrix $\Sigma$. }

First, suppose that the matrix $ \left[\begin{smallmatrix}
	A&B\\
	C&0
\end{smallmatrix}  \right]$ is PSD. Let   
$$
{\mathcal P}_1=\begin{bmatrix}
	0&0\\0&D
\end{bmatrix}, \quad\mathcal{P}_2=\begin{bmatrix}
	A&B\\
	C&0
\end{bmatrix},
\quad  \Sigma=	
\begin{bmatrix}
	\alpha Q_1&0\\0&\alpha Q_2
\end{bmatrix}.
$$
Hereafter, the corresponding iteration method is called  SPPS1 and the induced preconditioner is denoted by   ${\mathcal P}_{\SPPS1}$. So    
\begin{align*}
	{\mathcal P}_{\SPPS1}^{-1} & =( \Sigma + {\mathcal P}_2  )^{-1}\Sigma  ( \Sigma+ {\mathcal P}_1  )^{-1}\\
	&=
	\begin{bmatrix}
		A+\alpha Q_1&B\\C&\alpha Q_2
	\end{bmatrix}^{-1}
	\begin{bmatrix}
		\alpha Q_1&0\\0& \alpha Q_2
	\end{bmatrix}\begin{bmatrix}
		\alpha Q_1&0\\0&D+\alpha Q_2
	\end{bmatrix} ^{-1}
	\\&=
	\begin{bmatrix}
		I&0\\-\frac{1}{\alpha} Q_2^{-1}C&I
	\end{bmatrix}
	\begin{bmatrix}
		S_1^{-1}&0\\0&I
	\end{bmatrix}
	\begin{bmatrix}
		I&-B\\0&I
	\end{bmatrix}
	\begin{bmatrix}
		I&0\\0& (D+\alpha Q_2)^{-1}
	\end{bmatrix} ,
\end{align*}
where  $S_1=A+\alpha Q_1-\frac{1}{\alpha}B Q_2^{-1}C$.   Therefore, an algorithm for computing $ [y_1;y_2]= {\mathcal P}_{\SPPS1}^{-1} [x_1,x_2] $ 	can be written as the following algorithm.  		
\begin{algorithm}[H]
	\caption{Solving $ {\mathcal P}_{\SPPS1} [y_1;y_2]=[x_1,x_2] $. \label{alg-PPS1}}
	\begin{algorithmic}[1]
		\STATE  Solve $(D+\alpha Q_2)v_2=x_2 $;\\
		\STATE Set $ v_1=x_1-Bv_2 $; \\
		\STATE Solve $ (A+\alpha Q_1- B(\alpha Q_2)^{-1}C) y_1=v_1$;\\
		\STATE  Solve $ y_2=v_2- (\alpha Q_2)^{-1}Cy_1 $;  
	\end{algorithmic}
\end{algorithm}

For the second choice,	suppose that the matrix 
$ \left[\begin{smallmatrix}
	0&B\\
	C&D
\end{smallmatrix}  \right]$ is PSD.
we set    
$$
{\mathcal P}_1=\begin{bmatrix}
	A&0\\0&0
\end{bmatrix},\quad \mathcal{P}_2=\begin{bmatrix}
	0&B\\
	C&D
\end{bmatrix},
~~\Sigma=	
\begin{bmatrix}
	\alpha Q_1&0\\0&\alpha Q_2
\end{bmatrix}. 
$$
Henceforward, the corresponding iteration method is called SPPS2.  For the associated preconditioner ${\mathcal P}_{\SPPS2}$,  we have 
\begin{align*}
	{\mathcal P}_{\SPPS2}^{-1} & =( \Sigma + {\mathcal P}_2  )^{-1}\Sigma  ( \Sigma+ {\mathcal P}_1 )^{-1}\\
	&=\begin{bmatrix}
		\alpha Q_1&B\\C&D+\alpha Q_2
	\end{bmatrix} ^{-1}	\begin{bmatrix}
		\alpha Q_1&0\\0& \alpha Q_2
	\end{bmatrix}
	\begin{bmatrix}
		A+\alpha Q_1&0\\0&\alpha Q_2
	\end{bmatrix}^{-1}
	\\&=
	\begin{bmatrix}
		I&-(\alpha Q_1)^{-1}B\\0&I
	\end{bmatrix}
	\begin{bmatrix}
		I&0\\0& S_2^{-1}
	\end{bmatrix}
	\begin{bmatrix}
		I&0\\-C &I
	\end{bmatrix}	\begin{bmatrix}
		(A+\alpha Q_1)^{-1}&0\\0 &I
	\end{bmatrix} ,
\end{align*}
where  $S_2=D+\alpha Q_2-\frac{1}{\alpha}CQ_1^{-1}B $.  
Using the above factorization we can state the following algorithm  for computing $ [y_1;y_2]= {\mathcal P}_{\SPPS2}^{-1} [x_1,x_2] $.	 
\begin{algorithm}[H]
	\caption{Solving $ {\mathcal P}_{\SPPS2} [y_1;y_2]=[x_1,x_2] $. \label{alg-PPS2}}
	\begin{algorithmic}[1]
		\STATE  Solve $(A+  \alpha Q_1)v_1=x_1 $;\\
		\STATE Set $ v_2=x_2-Cv_1 $; \\
		\STATE Solve $ (D+  {\alpha} Q_2- C	(\alpha Q_1)^{-1}B ) y_2=v_2$;\\
		\STATE  Solve $ y_1=v_1- 	(\alpha Q_1)^{-1}B y_2 $;   
	\end{algorithmic}
\end{algorithm} 

 Considering the discussion given in   Section \ref{sec:sigma},   we may set 
$$\Sigma=\alpha_{*} Q,  ~~~\text{where}~  \quad 	\alpha_{*}=\frac{\norm{ {\mathcal A}}_F}{2\norm{Q}_F}, \quad \text{and} \quad Q\approx \frac{1}{2} ({\mathcal A}+{\mathcal A}^*) .$$ 
In Algorithms \ref{alg-PPS1} and \ref{alg-PPS2}, the matrices $ Q_1 $ and $ Q_2 $  should be HPD and the subsystems with them, as the coefficient matrix, should be solved easily.  Therefore,  if $A$ is PD, then we can choose $Q_1$ as an approximation of  $ \frac{1}{2}	(A+A^*)$ and if $A$ is PSD, then   we may take $ Q_1 $ as an approximation of  $ \frac{1}{2}(A+A^*)+\epsilon I$, where $\epsilon$ is a small positive number. In the same manner, the matrix $Q_2$ is chosen using the matrix $D$.  

In Algorithm \ref{alg-PPS1}, if the inverse of  $  Q_2$ is easily  computed,  then the computation   of the Schur complement involves simple matrix multiplication and subtraction operations. 
This can lead to efficient numerical computations and better performance compared to the  HSS, SS, TSS and BTSS methods. 
In practice, we solve the systems with  the Schur complement matrix   by an iterative method. So it is not necessary to form the Schur complement matrix explicitly.
Similarly, the matrix  $ Q_1 $ is selected in Algorithm~\ref{alg-PPS2}.

For example, in the cases that the matrices  $A$ and $C$ are, respectively, PD and PSD, then in Algorithm \ref{alg-PPS1},   we may consider $Q_1 = \frac{1}{2}(A+A^*)$ and $Q_2=\frac{1}{2} \diag(C+C^*)+\epsilon I$, and  in Algorithm \ref{alg-PPS2},    we can consider $Q_1 = \frac{1}{2} \diag(A+A^*)$ and $Q_2=\frac{1}{2}(C+C^*)+\epsilon I $. 

 {We present three numerical examples to show the effectiveness of our preconditioners. To do so, the numerical results of the  preconditioners SPPS1 and  SPPS2 are compared with those of the  TSS,   BTSS,  HSS,  SS, ETSS, EBTSS, EHSS, ESS,   and ILUT  (incomplete LU factorization) preconditioners  with dropping tolerance $droptol=0.1$. Since the  system  in Example \ref{test2} can also be solved using the modified HSS (MHSS) iteration method \cite{bai872010}, we also report numerical results of the FGMRES method in conjunction with the MHSS preconditioner for this example. Moreover, the parameter $\alpha$ is estimated 
	as  $\alpha=\sqrt{\gamma_{\min}\gamma_{\max}}$, where  $\gamma_{\min}$ and $\gamma_{\max}$ are  the extreme eigenvalues of  $ W $ (see \cite{bai872010}). With the aim of revealing  the superiority of the proposed preconditioners to the other preconditioners,
	the elapsed CPU time and the iterations  to achieve a solution of the system with a desired accuracy are compared.
} 

\begin{example}[\cite{bai872010}]\label{test2}
	Consider  the system of linear equations   of the form
	\begin{align} \label{sys2}
		(W + i T)x=b, \qquad  W,T \in \mathbb{R}^{p\times p},\qquad x,b\in \mathbb{C}^p,
	\end{align}
	where 
	$$ W= K+\frac{3-\sqrt{3}}{\tau}I,\quad \rm{and} \quad T=K+\frac{3+\sqrt{3} }{\tau} I, $$ 
	in which $ \tau $ is the time step-size and $ K $ is the five-point centered difference matrix approximating the negative Laplacian operator 
	$ L =-\Delta $ with homogeneous Dirichlet boundary conditions, 
	on a uniform mesh in the unit square $ [0, 1]  \times[0, 1] $ with the mesh-size
	$ h = \frac{1}{m+1} $.
	The matrix $ K \in \mathbb{R}^{p\times p}  $   possesses the tensor-product form $ K=I_m \otimes V_m + V_m \otimes I_m$, 
	with $ V_m = h^{-2}\text {tridiag}(-1, 2, -1) \in  \mathbb{R}^{m\times m} $.  Hence,  $ K$ is an $ p\times p$ block-tridiagonal   matrix, with $ p = m^2 $. 
	
	System \eqref{sys2} can be rewritten in the real form
	\begin{align}
		\mathcal{A} {u}=	\begin{bmatrix}
			W&-T\\T&W
		\end{bmatrix}
		\begin{bmatrix}  
			\Re(x)\\\Im(x)
		\end{bmatrix}=\begin{bmatrix}
			\Re(b)\\\Im(b)
		\end{bmatrix} =\textbf{b},
	\end{align}
	We consider this system for our numerical test.	The parameter  $\tau$ is set to be  $ \tau=h $ and  the $j$-th entry of the right-hand side vector  $  b $  
	is set to be
\begin{equation*} 
		b_j=\frac{(1- i) j}{\tau ( j + 1)^2}, \qquad  j = 1, 2, \ldots, p.
\end{equation*}
	In this example, we set $m=64$, $ 128 $, $256$ and  $512.$

\end{example}

\begin{example}\label{test5}
	We consider the Oseen problem
	\begin{align}
		\left\{\begin{aligned}\label{stoks}
			-v \Delta \mathbf{u}+\mathbf{w} \cdot \nabla \mathbf{u}+\nabla p=\mathbf{f}, & \text { in } \Omega, \\
			\nabla \cdot \mathbf{u}=0, & \text { in } \Omega,
		\end{aligned}\right.
	\end{align}
	with suitable boundary conditions on $\partial \Omega$, where $\Omega \subseteq \mathbb{R}^2$ is a bounded domain and $\mathbf{w}$ is a given divergence free field. The parameter $v>0$ is the viscosity, the vector field $\mathbf{u}$ stands for the velocity and $p$ represents the pressure. Here $\Delta, \nabla \cdot$  and $\nabla$ stand for the Laplace operator, the divergence operators and the gradient, respectively. The Oseen problem \eqref{stoks} is obtained from the linearization of the steady-state Navier-Stokes equation by the Picard iteration where the vector field $\mathbf{w}$ is the approximation of $\mathbf{u}$ from the previous Picard iteration. It is known that many discretization schemes for \eqref{stoks} will lead to a generalized saddle point problems of the form \eqref{sad}. We use the stabilized (the Stokes stabilization) $Q_1-P_0$ finite element method for the leaky lid driven cavity problems on uniform-grids on the unit square, with $v=0.001$. We use the IFISS software package developed by Elman et al. \cite{elman332007} to generate the linear systems corresponding to $32 \times 32,64 \times 64,128 \times 128$ and $256 \times 256$ grids.   
	In this example, the right-hand side vector $b$ is set to be $b=\mathcal{A} e$, where $e=[1,1, \ldots, 1]^T$.
\end{example}
\begin{example}\label{test6}
	In this example, we discuss the image restoration problems   $K X=F$, where $K \in \mathbb{R}^{p\times p} $ is a blurring matrix and $F=K X_{\text {org }}+\Phi$ is the observed image. Here, $X_{\text {org }} \in \mathbb{R}^{p\times q}$ and $\Phi \in \mathbb{R}^{p\times q}$ represent the original image and noise matrix, respectively. Consider solving the least squares problems with regularization
	\begin{align}\label{rest}
		\min _X\left\{\|K X-F\|_F^2+\mu\|G X\|_F^2\right\},
	\end{align}
	to restore the original image, where $\mu>0$ and $D$ is an auxiliary operator so that it reduces the effect of the noise vector and it is  often an approximation to a derivative 	operator. 
	It is not difficult to see that the problem  \eqref{rest} is equivalent to solve the saddle point problem 
\begin{equation*} 
	\begin{bmatrix}
		I_q \otimes I_p &I_q \otimes K\\-I_q \otimes K^*&\mu^2 I_q \otimes\left(G^* G\right)
	\end{bmatrix}\begin{bmatrix}
		y_1\\y_2
	\end{bmatrix}=\begin{bmatrix}
		\operatorname{vec}(F) \\0
	\end{bmatrix},
\end{equation*}
	where $\otimes$ is the Kronecker product $\operatorname{symbol}$ and $\operatorname{vec}(F)=\left[F_1 ; F_2 ; \cdots ; F_q\right]$, in which $F_i$ is the $i$-th column of $F$. In \eqref{rest}, let $K=\left(k_{i j}\right)$ be the $p \times p$  Toeplitz matrix with the entries 
	(see \cite{benzi272006})
	\begin{equation*}
		k_{i j}=\frac{1}{\sigma \sqrt{2 \pi}} \exp \left(\frac{-|i-j|^2}{2 \sigma^2}\right), \quad 1 \leqslant i, j \leq p,
	\end{equation*}
	and $G=\left(g_{i j}\right)$, $1\leq i,j\leq p$  where
	(see [50, pp 92-93])
	\begin{equation*}
		g_{ij}=\begin{cases}
			-2; &\text{if}~~~ i=j,\\
			1;&\text{if}~~ |i-j|=1 ~~\text{or} ~~|i-j|=p-1,\\
			0, &\text{otherwise.}
		\end{cases}
	\end{equation*}
	It is noted that if $\sigma$ is close to 2, then the matrix $K$ is a highly ill-conditioned  \cite{benzi272006}. 
	
	In this example, we set $\sigma=2,~ \mu=10^{-6}$ and noise $\Phi$  is set to be the Gaussian white noise with noise-to-signal ratios of 30dB. In addition, $X_{\text {org }}$ is the image  {\verb|cameraman.tft|}   of \textsc{Matlab}, which is of  $ 256\times 256 $ pixel image.

\end{example}

Generic properties of     matrices $\mathcal{A}$, $A$ and $D$    in \eqref{sad} for Examples \ref{test2}, \ref{test5} and \ref{test6},  have been presented in Table \ref{matrixpro}.
In Algorithm \ref{alg-PPS1}, it is necessary to solve a system  with $S_1=A+\alpha Q_1-\frac{1}{\alpha}B Q_2^{-1}C$ 
and in   Algorithm \ref{alg-PPS2}, a system with $S_2=D+\alpha Q_2-\frac{1}{\alpha}CQ_1^{-1}B $.  Therefore, if   the  non-zero entries of the matrix $ A $ in \eqref{sad} is less than that of $ D $, then we suggest to use the SPPS1 method with the Algorithm \ref{alg-PPS1}, otherwise   the SPPS2 method with the Algorithm \ref{alg-PPS2}. 
In Example \ref{test2}, we have   $ \nnz(A)=\nnz(D) $. So, there is no difference between using SPSS1 or  SPSS2 preconditioner. In Example \ref{test5}, since  $ \nnz(A)>\nnz(D) $, we suggest to use the SPSS2 preconditioner and in Example \ref{test6}, inasmuch as 
$ \nnz(A)<\nnz(D) $, we suggest to use the SPSS1 preconditioner.

We choose $\Sigma=\alpha_{*} Q$, where $\alpha_{*}$  is   defined  by \eqref{alpha2star} and $Q$ is  as an approximation of the matrix $ {\mathcal H}=\frac{1}{2} ({\mathcal A}+{\mathcal A}^*) $. 
With this in mind, let $H_A=\frac{1}{2}	(A+A^*)$ and  $H_D=\frac{1}{2}	(D+D^*)+\epsilon I$, where $\epsilon$ is equal to $0$, $10^{-2}$ and $10^{-4}$, respectively, in Examples \ref{test2}, \ref{test5} and \ref{test6}. 
For the SPPS1 and  SPPS2 preconditioners,  we set  $\Sigma$, respectively,    as
\begin{equation*}
	\Sigma=\alpha\begin{bmatrix} H_A&0\\
		0& \diag(H_D)
	\end{bmatrix}, \qquad \text{and } \qquad 	\Sigma=\alpha\begin{bmatrix} \diag(H_A)&0\\
		0& H_D
	\end{bmatrix}.
\end{equation*}  
For  {other  preconditioners}, we  utilize $\Sigma=\alpha Q$, where the matrix $Q$ is  of  the form  
\begin{equation*}
	\mathcal I ,   \quad
	\mathcal     D=
	\begin{bmatrix}
		\diag(H_A)&0\\
		0&\diag (H_D)
	\end{bmatrix} , \quad   \mathcal B=\begin{bmatrix}	H_A&0\\
		0&\diag (H_D)
	\end{bmatrix}.  
\end{equation*}
\begin{table}[t]
	\centering
	\begin{center}
		\caption{ Generic properties of     matrices $\mathcal{A}$, $A$ and  $D$    in \eqref{sad} for Examples \ref{test2}, \ref{test5} and \ref{test6}. \label{matrixpro}}
		\scriptsize
		\label{tb}
		{\footnotesize   
			\begin{tabular}{|c|c c|c|c|c|c|c|c|c|c|c|c|c|c|c|c|}
				\hline
				& &&$n$& cond($\mathcal{A}$)	&$ A $&$ D $&nnz( $A$ )& nnz( $D$ )\\\hline
				\multirow{4}{*}{Example \ref{test2}	}&
				\multirow{4}{*}{$ m$}&64&8192& 2.83e+2&\multirow{4}{*}{HPD}
				&\multirow{4}{*}{HPD}&20224&20224\\ 
				&&128&32768& 5.61e+2&  &&	81408&81408\\ 
				&&258&131072& 1.12e+3& &	&326656&326656\\ 
				&&512&524288& 2.23e+3& &&	1308672&1308672 \\ \hline
				\multirow{4}{*}{Example \ref{test5}}&\multirow{4}{*}{ grid}
				&32$ \times$ 32&3200&  9.88e+03&\multirow{4}{*}{PD} &\multirow{4}{*}{HPSD}&16818&3064 \\  
				&&64$ \times$ 64&12544&8.62e+04&&&70450&12280 \\  
				&&128 $ \times$128&49644& 1.38e+06 &&&288306&49144 \\  
				&&256$ \times$ 256&197632& 2.20e+07&&&1166386&196600\\\hline 
				\multirow{1}{*}{Example \ref{test6}}&\multirow{1}{*}{ pixel}
				&256$ \times$ 256&131072&6.84e+{14}&HPD&HPSD&65536&327680  \\\hline 
			\end{tabular}
		}
	\end{center}
\end{table}

\begin{table} 
	
	\begin{center}
		\caption{Numerical results  of the restarted FGMRES(30) for solving the preconditioned systems  for  Example \ref{test2}.}
		\scriptsize
		\label{test:2}
		{ \tiny
			\begin{tabular} { |c|c| c| c| c|c|c| c|c|c| c|c| }
				\hline
				$m$&splittig&SPPS1&SPPS2&TSS&BTSS, HSS&SS&ILUT& {MHSS}
				\\ 
				\hline \hline 
				\multirow[c]{3}{*}{\bfseries 64}&$\alpha_{*}$&0.75&$0.75$&$13463.63$&$13463.63$&$13463.63$&$-$&1859.88
				\\  
				&CPU&$0.25$&$0.15$&$1.44$&$1.01$&$1.42$&$0.29$&0.47
				\\
				&IT&$30$&$29$&$202$&$152$&$215$&$36$&144
				\\
				\hline \hline 
				\multirow[c]{3}{*}{\bfseries 128}&$\alpha_{*}$&0.75&$0.75$&$52827.66$&$52827.66$&$52827.66$&$-$&4942.59
			\\
				&CPU&$0.97$&$1.09$&$12.64$&$3.76$&$6.79$&$6.46$&2.6
		\\
				&IT&$42$&$40$&$408$&$279$&$420$&$57$&211
	\\
				\hline \hline 
				\multirow[c]{3}{*}{\bfseries 256}&$\alpha_{*}$&0.75&$0.75$&$209271.87$&$209271.87$&$209271.87$&$-$&13517.36
				\\
				&CPU&$8.61$&$8.92$&$126.46$&$28.31$&$52.55$&$143.15$&29.50
				\\
				&IT&$58$&$56$&$900$&$554$&$901$&$82$&425
				\\
				\hline \hline 
				\multirow[c]{3}{*}{\bfseries 512}&$\alpha_{*}$&0.75&$0.75$&$-$&$833024.81$&$-$&$-$&37568.92
				\\
				&CPU&$84.51$&$87.21$&$\ddag$&$310.55$&$\ddag$&$\ddag$&207.29
				\\
				&IT&$82$&$87$&$-$&$1215$&$-$&$-$&510
			\\
				\hline 
			\end{tabular}
		}\\
	\vspace{1mm}
	{ \tiny
		\begin{tabular} { |c|c| c|c|c| c|c|c|c|c|}
			\hline			
			$m$&splittig 
			& {ETSS-${\mathcal{D}}$}& {EBTSS-${\mathcal{D}}$, EHSS-${\mathcal{D}}$}& {ESS-${\mathcal{D}}$}& {ETSS-${\mathcal{B}}$}& {EBTSS-${\mathcal{B}}$, EHSS-${\mathcal{B}}$}& {ESS-${\mathcal{B}}$}\\ 
			\hline \hline 
			\multirow[c]{3}{*}{\bfseries 64}&$\alpha_{*}$ 
			&$0.79$&$0.79$&$0.79$&$0.75$&$0.75$&$0.75$ \\  
			&CPU  
			&$1.33$&$0.49$&$0.72$&$15.26$&$0.38$&$0.29$\\
			&IT 
			&$202$&$152$&$215$&$275$&$27$&$31$\\
			\hline \hline 
			\multirow[c]{3}{*}{\bfseries 128}&$\alpha_{*}$ 
			&$0.79$&$0.79$&$0.79$&$0.75$&$0.75$&$0.75$\\
			&CPU 
			&$12.08$&$2.79$&$5.85$&$112.35$&$2.94$&$1.84$\\
			&IT  
			&$408$&$279$&$420$&$748$&$38$&$43$\\
			\hline \hline 
			\multirow[c]{3}{*}{\bfseries 256}&$\alpha_{*}$ 
			&$0.79$&$0.79$&$0.79$&$-$&$0.75$&$0.75$\\
			&CPU 
			&$115.34$&$28.48$&$53.76$&$\ddag$&$15.95$&$16.32$\\
			&IT 
			&$900$&$556$&$901$&$-$&$60$&$59$\\
			\hline \hline 
			\multirow[c]{3}{*}{\bfseries 512}&$\alpha_{*}$ 
			&$-$&$0.79$&$-$&$-$&$0.75$&$0.75$\\
			&CPU 
			&$\ddag$&$310.04$&$\ddag$&$\ddag$&$325.62$&$169.53$\\
			&IT 
			&$-$&$1219$&$-$&$-$&$212$&$81$\\
			\hline 
		\end{tabular}
	}
\\
\vspace{5mm}
		\caption{Numerical results  of the restarted FGMRES(30) for solving the preconditioned systems  for Example  \ref{test5}.}
		
		\label{test:5}
		{ \tiny
			\begin{tabular}{ |c|c| c| c| c|c|c| c|c|}
				\hline 
				grid&splittig&SPPS1&SPPS2&
				TSS&BTSS&HSS&SS&ILUT
				 \\ 
				\hline \hline 
				\multirow[c]{3}{*}{\bfseries 32$\times$32}&$\alpha_{*}$&$0.51$&$0.51$&$0.15$&$0.15$&$0.15$&$0.15$&$-$ 
				\\  
				&CPU&$0.25$&$0.14$&$1.73$&$3.93$&3.72&$3.18$&$\ddag$
				\\  
				&IT&$43$&$41$&$1399$&$1625$&$1479$&$1701$&$-$
				\\  
				\hline \hline 
				\multirow[c]{3}{*}{\bfseries 64$\times$64}&$\alpha_{*}$&0.51&$0.51$&0.10&$0.10$&$0.10$&$0.10$&$-$
				\\  
				&CPU&$3.82$&$0.97$&$62.78$&$63.59$&65.48&$58.12$&$2.64$
				\\  
				&IT&$98$&$72$&$12199$&$12796$&12344&$12557$&$1631$
				\\  
				\hline \hline 
				\multirow[c]{3}{*}{\bfseries 128$\times$128}&$\alpha_{*}$&0.50&$0.50$&$-$&$-$&$-$&$-$&$-$
				\\  
				&CPU&$35.38$&$8.12$&$-$&$-$&$-$&$-$&$-$
				\\  
				&IT&$304$&$133$&$\dag$&$\dag$&$\dag$&$\dag$&$\dag$
				\\  
				\hline \hline 
				\multirow[c]{3}{*}{\bfseries 256$   { \times }$256}&$\alpha_{*}$&0.50&$0.50$&$-$&$-$&$-$&$-$&$-$
				\\  
				&CPU&$456.63$&$192.62$&$\ddag$&$\ddag$&$\ddag$&$\ddag$&$\ddag$
				\\  
				&IT&$642$&$862$&$-$&$-$&$-$&$-$&$-$
				\\  
				\hline 
			\end{tabular}
		}
	\\\vspace{1mm}
%
	{ \tiny
		\begin{tabular}{ |c|c| c| c| c|c|c| c|c|c| c|c|c| c|c|c|c|c|c|}
			\hline 
			grid&splittig& 
			{ETSS-${\mathcal{D}}$}& {EBTSS-${\mathcal{D}}$}& {EHSS-${\mathcal{D}}$}& {ESS-${\mathcal{D}}$}&
			{ETSS-${\mathcal{B}}$}& {EBTSS-${\mathcal{B}}$}& {EHSS-${\mathcal{B}}$}& {ESS-${\mathcal{B}}$} \\ 
			\hline \hline 
			\multirow[c]{3}{*}{\bfseries 32$\times$32}&$\alpha_{*}$ 
			&$0.51$&$0.51$&$0.51$&$0.51$&$0.51$&$-$&$-$&$0.51$\\  
			&CPU 
			&$1.63$&$0.79$&0.71&$0.44$&$37.89$&$-$&$-$&$0.35$\\  
			&IT 
			&$81$&$181$&220&$87$&$3419$&$\dag$&$\dag$&$27$\\  
			\hline \hline 
			\multirow[c]{3}{*}{\bfseries 64$\times$64}&$\alpha_{*}$ 
			&$0.51$&$0.51$&0.51&$0.51$&$-$&$-$&$-$&$0.51$\\  
			&CPU 
			&$14.61$&$10.16$&15.32&$3.85$&$\ddag$&$-$&$-$&$4.35$\\  
			&IT 
			&$257$&$919$&1200&$209$&$-$&$\dag$&$\dag$&$74$\\  
			\hline \hline 
			\multirow[c]{3}{*}{\bfseries 128$\times$128}&$\alpha_{*}$ 
			&$0.50$&$-$&$-$&$0.50$&$-$&$-$&$-$&$0.50$\\  
			&CPU 
			&$138.38$&$-$&$-$&$32.00$&$\ddag$&$\ddag$&$\ddag$&$48.31$\\  
			&IT 
			&$776$&$\dag$&$\dag$&$640$&$-$&$-$&$-$&$271$\\  
			\hline \hline 
			\multirow[c]{3}{*}{\bfseries 256$   { \times }$256}&$\alpha_{*}$ 
			&$-$&$-$&$-$&$-$&$-$&$-$&$-$&$-$\\  
			&CPU 
			&$\ddag$&$\ddag$&$\ddag$&$\ddag$&$\ddag$&$\ddag$&$\ddag$&$\ddag$\\  
			&IT 
			&$-$&$-$&$-$&$-$&$-$&$-$&$-$&$-$\\  
			\hline 
		\end{tabular}
	}
\\
\vspace{5mm}
		\caption{Numerical results  of the restarted FGMRES(30) for solving the preconditioned systems  for the Example \ref{test6}.}
		
		\label{test:6}
		{ \tiny  
			\begin{tabular}{|c|c| c| c| c|c|c| c|c|c| c|c|c| c|c|c|}
				\hline
						pixel&splittig&SPPS1&SPPS2&
				TSS&BTSS{, HSS}&SS& ILUT\\ 
				\hline \hline 
				\multirow[c]{5}{*}{\bfseries 256$ \times$256 }&$\alpha_{*}$&0.57&$0.57$&$0.40$&$0.40$&$0.40$&$-$
				\\  
				&CPU&$69.48$&$70.11$&$297.85$&$324.44$&$379.65$&$\ddag$
				\\  
				&IT&$79$&$59$&$840$&$883$&$1057$&$-$
				\\  
				&SNR&$24.77$&$24.76$&$24.77$&$24.77$&$24.77$&$-$
				\\  
				&PSNR&$30.35$&$30.34$&$30.35$&$30.35$&$30.35$&$-$
				\\  
				\hline 
			\end{tabular}
		}
	\\\vspace{1mm}
%
	{ \tiny  
		\begin{tabular}{|c|c| c| c| c|c|c| c|c|c| c|c|c| c|c|c|}
			\hline
			pixel&splittig& {ETSS-${\mathcal{D}}$}& {EBTSS-${\mathcal{D}}$, EHSS-${\mathcal{D}}$}& {ESS-${\mathcal{D}}$}& {ETSS-${\mathcal{B}}$}& {EBTSS-${\mathcal{B}}$, EHSS-${\mathcal{B}}$}& {ESS-${\mathcal{B}}$}\\ 
			\hline \hline 
			\multirow[c]{5}{*}{\bfseries 256$ \times$256 }&$\alpha_{*}$ 
			&$-$&$-$&$0.57$&$-$&$-$&$0.57$\\  
			&CPU 
			&$-$&$\ddag$&$89.51$&$-$&$\ddag$&$88.75$\\  
			&IT 
			&$\S$&$-$&$82$&$\S$&$-$&$82$\\  
			&SNR 
			&$-$&$-$&$24.77$&$-$&$-$&$24.77$\\  
			&PSNR 
			&$-$&$-$&$30.35$&$-$&$-$&$30.35$\\  
			\hline 
		\end{tabular}
	}
	\end{center}
\end{table}

\begin{figure}[!t] 
	\centering
	\includegraphics[width=0.84\linewidth]{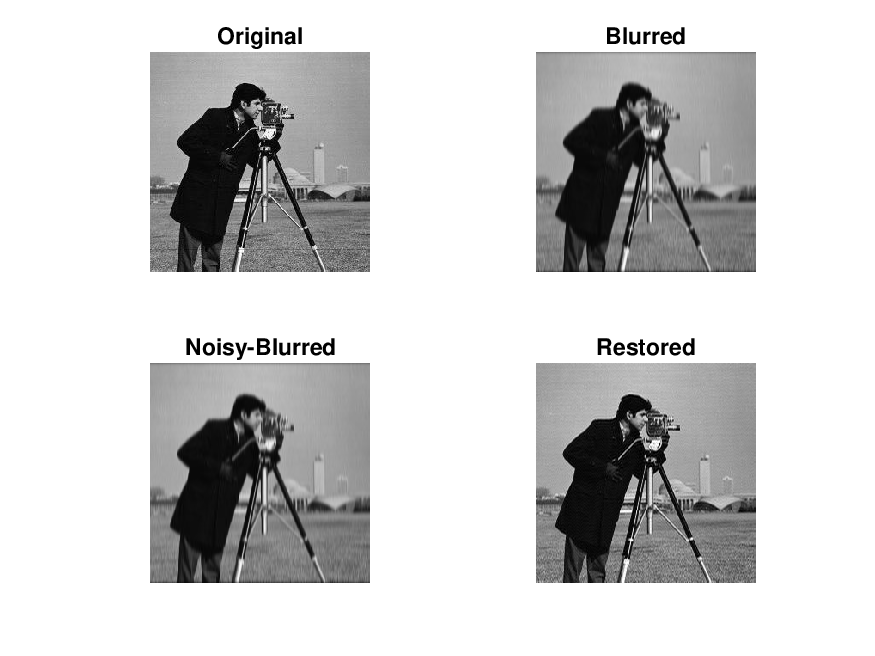}
	\caption{Original, blurred, noisy-blurred and restored images of Cameraman. \label{fig-cameraman}} 
\end{figure}

We utilize the restarted version of the flexible GMRES(30) (FGMRES(30))   method alongside the aforementioned preconditioners to solve the given examples. The outer iteration is terminated once the Euclidean norm of the residual is decreased by a factor of $10^7$. For the subsystems with HPD coefficient matrices, we utilize the CG method, while those with PD matrices are solved using the restarted GMRES(10).  For the inner systems, the iterations are stopped when the residual 2-norm is reduced by a factor of 
$10^1$, and a maximum of 50 iterations is allowed. All runs are carried out in \textsc{Matlab} 2017a on an Intel core(TM) i7-8550U (1.8 GHz) 16G RAM Windows 10 system.

Tables \ref{test:2},  \ref{test:5}  and  \ref{test:6}  display  the numerical results along with the parameter $\alpha_{*}$, for  Examples \ref{test2}, \ref{test5} and \ref{test6}.  
 {In the tables,  the methods with the shift matrices  $\Sigma=\alpha \mathcal{D}$ and 
	$\Sigma=\alpha \mathcal{B}$ are  abbreviated as 	``Method's name"- $\mathcal{D}$ and  ``Method's name"- $\mathcal{B}$, respectively.}
These tables also present the elapsed CPU time (CPU) and the number of iterations (IT) for the convergence.  In all the tables, the notation $\dag$ means that the method has 	not converged in at most 15000 iterations and a $\ddag$ shows that method has not converged in at most 500s. Moreover, notation $\S$ shows that the  method fails to converge. 	In Table \ref{test:6},   the signal-to-noise ratio (SNR)  and peak signal-to-noise ratio (PSNR)  are presented, where 
\begin{equation*}
	{\rm SNR}=20\log_{10}  \left(\frac{\norm{X_{\text {org }}}}{\norm{X_{\text {org }}-X }}\right) , ~~~~~~
	{\rm PSNR}=20\log_{10} \left(\frac{   255\sqrt{pq} }{\norm{X_{\text {org }} -X}}\right).
\end{equation*}
where, $p=q=256$.
The original, blurred, noisy-blurred and the restored images  of   {\verb|cameraman|}  have been shown in Fig. \ref{fig-cameraman}, when the  SPPS1 preconditioner is employed.  
When the picture is restored successfully, there is not any   
significant difference between restored images, however from the CPU time point of view, the SPPS1 preconditioner outperforms the others. 

 	The tables indicate that as the matrix $Q $ approaches the Hermitian part of the matrix $\mathcal A$, the speed of convergence typically increases. As the numerical results show  the {SPPS1} and {SPPS2} preconditioners outperform the others. Notably, in many cases, the restarted FGMRES(30) method with the TSS,  BTSS,  HSS,  SS, ETSS, EBTSS, EHSS, ESS,   and ILUT  preconditioners failed to provide any solution, whereas the methods using the {SPPS1} and {SPPS2} preconditioners consistently succeeded.

\section{Conclusion}\label{sec:conc}
In this article, we have introduced the PPS method  for solving the system of linear equations with nonsingular  PSD coefficient matrices. The PPS method encompasses and extends existing methods such as SS, HSS, NSS, and PSS. 
In this approach, we split the coefficient matrix $\mathcal{A}$ into $\mathcal{A} = {\mathcal P}_1 + {\mathcal P}_2$, where both ${\mathcal P}_1$ and ${\mathcal P}_2$ are PSD matrices. Based on this splitting, we have introduced the preconditioner ${\mathcal P}_{\PPS} = (\Sigma + {\mathcal P}_1)\Sigma^{-1}(\Sigma + {\mathcal P}_2)$, where $\Sigma$ is an HPD matrix.
We have studied the convergence  properties of this method, providing necessary and sufficient conditions for its convergence.  It has been observed that when certain conditions are met, such as the positive definiteness or skew-Hermitian  or Hermitian of one of the ${\mathcal P}_1$ or ${\mathcal P}_2$, or    $ \nul({\mathcal P}_1)\cup \nul({\mathcal P}_2+ {\mathcal P}_2^*)=\mathbb{C} ^n$ or $ {\mathcal P}_1 \Sigma ^{-1}  {\mathcal P}_2 = {\mathcal P}_2 \Sigma ^{-1} {\mathcal P}_1 $, then this method is convergent. But in general, we have $\rho(\Gamma_{\PPS})\leq 1$. 
Nevertheless, the preconditioner ${\mathcal P}_{\PPS}$
can be effectively applied for solving linear systems with nonsingular PSD coefficient matrices. We also have demonstrated that the   matrix $ \Sigma $ can any arbitrary HPD matrix, significantly influencing the spectral radius of the iteration matrix and the convergence speed of the preconditioning method. When  $ \Sigma=\alpha Q $, a  specific value for the parameter  $\alpha$ was proposed. Additionally, we provided guidance on the selection of the matrices $ {\mathcal P}_1 $ and $ {\mathcal P}_2 $.  Furthermore, we discussed two special cases of the PPS method, showing that these variants are superior to methods like SS, HSS, TSS, and BTSS in terms of performance.

Finally, it is worth noting that one may generalize the PPS method by selecting two different HPD shift matrices  $\Sigma_1$ and $\Sigma_2$ in two half-steps of the method. In this case,  if the matrix ${\mathcal P}_2$ is skew-Hermitian and  ${\mathcal P}_1$ is an HPD matrix, then the method reduces to GPHSS presented in \cite{yang2162010} if we set $\Sigma_1=\alpha \mathcal P$  and $\Sigma_2=\beta \mathcal P$,  and to the 
GHSS method given in \cite{li2362012} if we choose $\Sigma_1=\alpha \mathcal I$ and $\Sigma_2=\beta \mathcal I$.
However, the convergence analysis of the method in the general form  would be  complicated and currently is under consideration. 

\section*{Acknowledgements}
The authors would like to thank the anonymous referee for the helpful comments and suggestions.	The authors would also like to thank Prof. Fatemeh Panjeh Ali Beik (Vali-e-Asr University of Rafsanjan) for careful proofreading of the paper. 
This work is supported by the Iran National Science Foundation (INSF) under Grant 	No. 4006300.



\end{document}